\documentclass[11pt,a4paper,notitlepage,final,onecolumn,oneside]{article}
\usepackage{exscale}
\usepackage{latexsym}
\usepackage{calc}
\usepackage{amssymb}
\usepackage{url}

\newcommand{\proof}{{\bfseries Proof:\quad}}

\begin{document}

\newtheorem{lemma}{Lemma}
\newtheorem{proposition}[lemma]{Proposition}
\newtheorem{theorem}[lemma]{Theorem}
\newtheorem{definition}[lemma]{Definition}
\newtheorem{hypothesis}[lemma]{Hypothesis}
\newtheorem{conjecture}[lemma]{Conjecture}
\newtheorem{remark}[lemma]{Remark}
\newtheorem{example}[lemma]{Example}
\newtheorem{property}[lemma]{Property}
\newtheorem{corollary}[lemma]{Corollary}
\newtheorem{algorithm}[lemma]{Algorithm}


\title{On the polycirculant conjecture}
\author{Aleksandr Golubchik\\
        \small (Osnabrueck, Germany)}

\maketitle

\begin{abstract}
In the paper the foundation of the $k$-orbit theory is developed. The
theory opens a new simple way to the investigation of groups and
multidimensional symmetries.

The relations between combinatorial symmetry properties of a $k$-orbit and
its automorphism group are found. It is found the local property of a
$k$-orbit. The difference between 2-closed group and $m$-closed group for
$m>2$ is discovered. It is explained the specific property of Petersen
graph automorphism group $n$-orbit. It is shown that any non-trivial
primitive group contains a transitive imprimitive subgroup and as a result
it is proved that the automorphism group of a vertex transitive graph
(2-closed group) contains a regular element (polycirculant conjecture).

Using methods of the $k$-orbit theory, it is considered different
possibilities of permutation representation of a finite group and shown
that the most informative, relative to describing of the structure of a
finite group, is the permutation representation of the lowest degree.
Using this representation it is obtained a simple proof of the W. Feit,
J.G. Thompson theorem: Solvability of groups of odd order. It is described
the enough simple structure of lowest degree representation of finite
groups and found a way to constructing of the simple full invariant of a
finite group.

To the end, using methods of $k$-orbit theory, it is obtained one of
possible polynomial solutions of the graph isomorphism problem.

\end{abstract}


\section{Introduction}
A permutation group $G$ on a $n$-element set $V$ is called regular, if its
every stabilizer (a subgroup that fixes some element $v\in V$) is trivial.
Every permutation of the regular group can be decomposed into cycles of
the same length.

A permutation group, containing a regular subgroup, is called semiregular.

Let $G$ be a permutation group on a $n$-element set $V$, $V^k$ be
Cartesian power of $V$ and $V^{(k)}\subset V^k$ be the non-diagonal part of
$V^k$, i.e.  every $k$-tuple of $V^{(k)}$ has $k$ different values of its
coordinates. The action of $G$ on $V$ forms the partition of $V^{(k)}$ on
classes of $k$-tuples related with $G$.  This partition is called the
system of $k$-orbits of $G$ on $V^{(k)}$ and we write it as $Orb_k(G)$. If
$\langle v_1\ldots v_k\rangle\in V^{(k)}$, then $G\langle v_1\ldots
v_k\rangle\equiv \{\langle gv_1\ldots gv_k\rangle : g\in G\}$ is a
$k$-orbit from $Orb_k(G)$.

For considered tasks it is of interest a maximal subgroup of
$Aut(Orb_k(G))$ that maintains $k$-orbits from $Orb_k(G)$. We shall denote
this subgroup as $aut(Orb_k(G))$. Thus $aut(Orb_k(G))\equiv\cap_{X_k\in
Orb_k(G)} Aut(X_k)$ and $G\leq aut(Orb_k(G))$.

\begin{definition}
We call a permutation group $G$ a \emph{$k$-defined} group, if
$G=aut(Orb_k(G))$.

\end{definition}

There is the obvious property of a $k$-defined group:

\begin{proposition}
If a group $G$ is $k$-defined, then it is $(k+1)$-defined.

\end{proposition}
\proof
If a group $G$ is $k$-defined, then, on the one hand,
$aut(Orb_{k+1}(G))<aut(Orb_k(G))=G$ and, on the other hand,
$G<aut(Orb_{k+1}(G)$. $\Box$\bigskip

The $k$-defined group is called $k$-closed if it is not $(k-1)$-defined.

P. Cameron \cite{Cameron} has described the conjecture of M. Klin, that
every $2$-closed transitive group is semiregular, and the similar
polycirculant conjecture of D. Maru$\breve{\rm s}$i$\breve{\rm c}$, that
every vertex-transitive finite graph has a regular automorphism.

We shall prove these conjectures in the next reformulation:

\begin{definition}\label{tr.pr}
We shall emphasize in the conventional definition of primitivity and
imprimitivity of permutation groups the case of cyclic group of a prime
order. We shall say that these groups are \emph{trivial primitive} and
\emph{trivial imprimitive}. The reason of such consideration will be clear
below.

\end{definition}

\begin{theorem}\label{2cl.impr->reg}
The $2$-closure of a transitive, imprimitive permutation group contains a
regular element.

\end{theorem}

\begin{lemma}\label{tr.impr<tr.pr}
A primitive permutation group contains a transitive, imprimitive subgroup.

\end{lemma}

In order to prove these statements we shall study symmetry properties of
$k$-orbits.\footnote{The objects, that generalize symmetry properties of
$k$-orbits, were applied by author for the polynomial solution of graph
isomorphism problem. The part of such investigations is used in
\url{http://arXiv.org/find/math/1/AND+au:+Golubchik_Aleksandr+ti:+AND+polynomial+algorithm/0/1/0/2002/0/1}}


\section{$k$-Orbits}
A $k$-orbit $X_k$ is a set of $k$-tuples with property
$X_k=Aut(X_k)\alpha_k$ for a $k$-tuple $\alpha_k\in X_k$. Such $k$-sets we
shall call \emph{automorphic} $k$-sets.

All, what is written below, can become easier for understanding, if to
represent a $k$-orbit as a matrix, whose lines are $k$-tuples and columns
are values of coordinates of $k$-tuples. A $k$-orbit can be represented by
various matrices that differ by lines permutation. Various orders of lines
in matrices demonstrate various symmetry properties of $k$-orbit. For
example $3$-orbit of symmetric group $S_3$ we can represent as

$$
\begin{array}{lcr}
\begin{array}{||ccc||}
\hline
1 & 2 & 3 \\
2 & 3 & 1 \\
3 & 1 & 2 \\
\hline
1 & 3 & 2 \\
3 & 2 & 1 \\
2 & 1 & 3 \\
\hline
\end{array} &
\mbox{ or } &
\begin{array}{||cc|c||}
\hline
1 & 2 & 3 \\
2 & 1 & 3 \\
\hline
1 & 3 & 2 \\
3 & 1 & 2 \\
\hline
2 & 3 & 1 \\
3 & 2 & 1 \\
\hline
\end{array}\,.
\end{array}
$$

In order to indicate number of $k$-tuples in a $k$-orbit $X_k$ of power
$l$ we shall call it $(l,k)$-orbit or write as $X_{lk}$.

$k$-Orbits have the next general number property:

\begin{proposition}\label{homo}
Let $X_k$ be a $k$-orbit and $l\leq k$, then all $l$-tuples with the same
coordinates $i_1,\ldots i_l\in [1,k]$ form a homogeneous multiset (i.e.
all $l$-tuples in this multiset have the same multiplicity).

\end{proposition}
\proof
Let two $k$-subsets $Y_k(u_1,u_2,\ldots,u_l),Z_k(v_1,v_2,\ldots,v_l)
\subset X_k$ consist of all $k$-tuples, that have their $l$ coordinates
$i_1,\ldots i_l$ equal to $u_1,u_2,\ldots,u_l$ and $v_1,v_2,\ldots,v_l$
correspondingly. Let $g\in Aut(X_k)$ be such permutation that $\langle
v_1\ldots v_l\rangle =\langle gu_1\ldots gu_l\rangle$, then $Z_k=gY_k$.
$\Box$\bigskip

Following constructions simplify the study of $k$-orbits. We call a
$k$-orbit as a \emph{cyclic $k$-orbit} or simply a \emph{$k$-cycle}, if it
is generated by single permutation. A $k$-cycle, that consists of $l$
$k$-tuples, we write as $(l,k)$-cycle. The order of a generating $k$-cycle
permutation can differ from the number of $k$-tuples in the $k$-cycle. The
structure of $k$-cycles is enough simple and can be represented with four
structure elements:

\begin{example}\label{simple-struc}
$$
\begin{array}{||cc|cc||}
\hline
1 & 2 & 3 & 4 \\
2 & 1 & 4 & 3 \\
\hline
\end{array}\,, \
\begin{array}{||cc|cc||}
\hline
1 & 2 & 3 & 4 \\
2 & 1 & 3 & 4 \\
\hline
\end{array}\,, \
\begin{array}{||cc|cc||}
\hline
1 & 2 & 3 & 4 \\
2 & 1 & 5 & 6 \\
\hline
\end{array}\,, \
\begin{array}{||cccc|cc||}
\hline
1 & 2 & 3 & 4 & 5 & 6 \\
2 & 3 & 4 & 1 & 6 & 5 \\
3 & 4 & 1 & 2 & 5 & 6 \\
4 & 1 & 2 & 3 & 6 & 5 \\
\hline
\end{array}\,.
$$

\end{example}

The first example shows the $(2,4)$-cycle that is a \emph{concatenation} of
two $(2,2)$-cycles. Such $(k,k)$-cycle, that is a $k$-orbit of a cycle of
length $k$, we shall call a \emph{$k$-rcycle}. This term is an
abbreviation from a ``right cycle'' and indicates on invariance of such
$(k,k)$-cycle relative to cyclic permutation of not only its $k$-tuples
but also the coordinates of $k$-tuples, or on invariance of the
$(k,k)$-cycle relative to not only the left but also the right action of
permutation (s. below).

The second is the $(2,4)$-cycle with fix-points. It is represented by the
concatenation of the $2$-rcycle and the trivial $(2,2)$-multiorbit,
consisting of the single $2$-tuple. Such $k$-multiorbit we shall call a
\emph{$k$-multituple} or $(l,k)$-multituple.

The third example is the concatenation of the $2$-rcycle and a $2$-orbit
that consists of two $2$-tuples with not intersected values of
coordinates. This kind of $(l,k)$-orbit we shall call
\emph{$S_l^k$-orbit}. It designates that this $(l,k)$-orbit consists of
$lk$ elements of $V$ and its automorphism group is subdirectproduct of
symmetric groups $S_l(B_i)$, where $B_i\subset V$, $i\in [1,k]$, $|B_i|=l$
and $B_i\cap B_j=\emptyset$ for $i\neq j$. From this definition follows
that any $(l,1)$-orbit ($l$-element set) is $S_l^1$-orbit.

The fourth example shows the possible structure of a $k$-cycle whose
length is not prime. It is seen that the fourth case can be represented
through first three cases. So these three cases are fundamental for
constructing of any $k$-orbit of any finite group.

One of our tasks is the study of a permutation action on $k$-orbits.
Indeed, there exist different possibilities of the permutation action on
$k$-orbits, which are arisen from their different symmetry properties.

We shall start with consideration of permutation actions on a $n$-orbit
$X_n$ of a group $G$ of the degree $n$.


\subsection{The actions of permutations on $n$-sets}

A $n$-orbit $X_n$ of a group $G$ is a set of all $n$-tuples, any pair of
which defines a permutation from $G$. So we can represent any $n$-tuple
$\alpha_n=\langle u_1\ldots u_n\rangle$ as a permutation
$$
g_{\alpha_n}=
\left( \begin{array}{ccc}
v_1 & \ldots & v_n \\
u_1 & \dots & u_n
\end{array}
\right),
$$
where $n$-tuple $\langle v_1\ldots v_n\rangle$ is related to the unit of
$G$ and will be called as the \emph{initial} $n$-tuple. Of course, any
$n$-tuple from $X_n$ can be chosen as the initial. The property of the
initial $n$-tuple is the equality of number value to order value of
each its coordinate. Here is accepted that sets of number values and order
values of coordinates are equal and for ordering of coordinates it is
determinate (any time if it is necessary) certain linear order on this
set. The next example shows two different orders of coordinates of the
same $2$-orbit:

$$
\begin{array}{lcr}
\begin{array}{||cc||}
\hline
1 & 2 \\
\hline
1 & 2 \\
2 & 1 \\
\hline
\end{array}\,,
&&
\begin{array}{||cc||}
\hline
2 & 1 \\
\hline
1 & 2 \\
2 & 1 \\
\hline
\end{array}\,.
\end{array}
$$

In first case the initial $2$-tuple is $\langle 12\rangle$, in the second
it is $\langle 21\rangle$.

Further we shall take the next rule for permutation multiplication:
$$
\left(
\begin{array}{ccc}
v_1 & \dots  & v_n \\
u_1 & \dots  & u_n
\end{array}
\right)
\left(
\begin{array}{ccc}
v_1 & \dots  & v_n \\
w_1 & \dots  & w_n
\end{array}
\right)
=
\left(
\begin{array}{ccc}
w_1 & \dots  & w_n \\
x_1 & \dots  & x_n
\end{array}
\right)
\left(
\begin{array}{ccc}
v_1 & \dots  & v_n \\
w_1 & \dots  & w_n
\end{array}
\right).
$$

From this rule follows that the left action of the permutation

$$
\left(
\begin{array}{ccc}
v_1 & \dots  & v_n \\
u_1 & \dots  & u_n
\end{array}
\right)
$$
on the $n$-tuple $\alpha_n=\langle w_1\dots w_n\rangle$ gives the $n$-tuple
$\beta_n=\langle x_1\dots x_n\rangle$ that can be considered as:

\begin{enumerate}
\item
the changing of (number) values of coordinates of the $n$-tuple $\alpha_n$;

\item
the mapping of $n$-tuple $\alpha_n$ coordinate-wise on a $n$-tuple
$\beta_n$.

\end{enumerate}

The right action of the permutation

$$
\left(
\begin{array}{ccc}
v_1 & \ldots & v_n \\
w_1 & \dots  & w_n
\end{array}
\right)
$$

on the $n$-tuple $\alpha_n=\langle u_1\dots u_n\rangle$ gives the $n$-tuple
$\beta_n=\langle x_1\dots x_n\rangle$ that can be interpreted as:

\begin{enumerate}
\item
the permutation of coordinates of the $n$-tuple $\alpha_n$;

\item
the mapping of $n$-tuple $\alpha_n$ coordinate-wise on a $n$-tuple
$\beta_n$.

\end{enumerate}

We shall choose every time such interpretation of permutation action that
will be more suitable.

If a $n$-orbit $X_n$ contains a $n$-tuple $\alpha_n=\langle v_1\ldots
v_n\rangle$, then $X_n=\{\langle gv_1\ldots gv_n\rangle : g\in G\}$ and
$|X_n|=|G|$.  Here we have used the first method of permutation action on
$n$-tuple, namely: a permutation $g$ changes values of coordinates of
$\alpha_n$ or acts on the permutation $g_{\alpha_n}$ from left
($gg_{\alpha_n}$). We shall say also that a permutation $g$ acts from left
on the $n$-tuple $\alpha_n$ and write this action as $g\alpha_n$.

The second method gives $X_n=\{\langle v_{g1}\ldots v_{gn}\rangle: g\in
G\}$. It is an action of a permutation $g$ on the order of coordinates of
$n$-tuple $\alpha_n$ or the action $g_{\alpha_n}g$ of $g$ on
$g_{\alpha_n}$ from right. We shall say in this case that $g$ acts on the
$n$-tuple $\alpha_n$ from right and write this action as $\alpha_ng$.

\paragraph{$n$-Orbits of left, right cosets of a subgroup.}

Let $A$ be a subgroup of $G$, $Y_n$ be a $n$-orbit of $A$ and $g\in G$,
then $gY_n$ is a subset of $X_n$ that represents permutations from a
left coset of $A$ in $G$ and a subset $Y_ng$ represents permutations from
a right coset of $A$ in $G$.

\begin{definition}\label{GA,AG}
For convenience we shall write the sets of left $G\backslash A$ and right
$G/A$ cosets of $A$ in $G$ as $GA\equiv \{gA:  g\in G\}$ and $AG\equiv
\{Ag: g\in G\}$.

The defined notation can be easily distinguished from the group
multiplication, because the product of a group with its subgroup is always
trivial. The same reasoning will be used in the like formulas by the
action of a group on $k$-orbits.

Corresponding to this remark we write $GY_n=\{gY_n: g\in G\}$,
$Y_nG=\{Y_ng: g\in G\}$, $AX_n=\{A\alpha_n: \alpha_n\in X_n\}$ and
$X_nA=\{\alpha_nA: \alpha_n\in X_n\}$.

\end{definition}\bigskip

From this definition and notions of left and right cosets of a subgroup we
obtain:

\begin{lemma}\label{GY_n=X_nA}
Let $n$-orbit $Y_n$ of a subgroup $A<G$ contains initial $n$-tuple, then
partitions of $X_n$ on $n$-subsets of left, right cosets of $A$ in $G$
are  $L_n=GY_n=X_nA$ and $R_n=Y_nG=AX_n$.

\end{lemma}\bigskip

Let $Y_n$ be a $n$-orbit of a subgroup $A$ of a group $G$ and $g\in G$.

\begin{lemma}\label{Y_ng}
The $n$-subset $Y_ng$ is also a $n$-orbit of the subgroup $A$.

\end{lemma}
\proof
The right action of a permutation on a $n$-orbit changes the order of
coordinates of $n$-tuples. The permutation that is defined by any pair of
$n$-tuples does not depend on order of coordinates, hence every
permutation that is defined by any pair of $n$-tuples from $Y_ng$ belongs
to $A$. $\Box$

\begin{lemma}\label{gY_n}
The $n$-subset $gY_n$ is a $n$-orbit of the conjugate to $A$ subgroup
$B=gAg^{-1}$.

\end{lemma}
\proof
The $n$-subsets $gY_n$ and $gY_ng^{-1}$ define, as in lemma \ref{Y_ng},
the same sets of permutations from $G$. But the set of $n$-tuples
$gY_ng^{-1}$ is by definition equivalent to the set of permutations
$B=gAg^{-1}$. $\Box$

\begin{proposition}\label{Yn'=gYn''f-1}
$n$-Subsets of left and right cosets of a subgroup $A<G$ are connected
with elements of $G$.

\end{proposition}
\proof Let $Y_n$ be a $n$-orbit of $A$, $Y_n'$ be the $n$-subset of
a left coset of $A$ and $Y_n''$ be the $n$-subset of a right coset of $A$,
then there exist permutations $g,f\in G$ so that $Y_n'=gY_n$ and
$Y_n''=Y_nf$. Hence $Y_n'=gY_n''f^{-1}$. $\Box$\bigskip

We shall say further ``$n$-orbit of coset'' instead of ``$n$-subset of
coset'', in order to show that this $n$-subset is a $n$-orbit. It will be
referred also to a $k$-subset, if it is a $k$-orbit.

\begin{proposition}\label{H->Yn=gYng-1}
Let $H$ be a normal subgroup of $G$, then sets of $n$-orbits of left and
right cosets of $H$ are equal and if to choose an arbitrary $n$-tuple from
$n$-orbit $Y_n$ of an arbitrary coset of $H$ as the initial, then
$Y_n=gY_ng^{-1}$ for any permutation $g\in G$.

\end{proposition}
\proof
The sets of $n$-orbits of left and right cosets of $H$ are equal, because
the sets of left and right cosets of $H$ are equal.

$Y_n=gY_ng^{-1}$ for every permutation $g\in G$, because the choice of the
initial $n$-tuple determines an equivalence between $Y_n$ and $H$
relative to the defined above action of $G$ on its $n$-orbit.
$\Box$\smallskip

\begin{proposition}\label{Yn=gYng-1->N}
Let $A$ be a subgroup of a group $G$, $Y_n$ be a $n$-orbit of $A$ and
$Y_n'$, $Y_n''$ be $n$-orbits of left and right cosets of $A$
correspondingly. Let $Y_n'=Y_n''\neq Y_n$, then $A$ has non-trivial
normalizer $N_G(A)$.

\end{proposition}
\proof
There exist permutations $g,f\in G\setminus A$ so that
$Y_n'=gY_n=Y_nf=Y_n''$. From this equality and lemma~\ref{Y_ng} follows
that $gY_ng^{-1}$ is a $n$-orbit of $A$. As $Y_n$ contains the initial
$n$-tuple, then $gY_ng^{-1}\cap Y_n\neq \emptyset$ and hence
$gY_ng^{-1}=Y_n$. So $gAg^{-1}=A$. $\Box$

\begin{proposition}\label{g*an*g-1,f*Yn*f-1}
Let $A$ be a subgroup of $G$, $g,f\in G$, $gag^{-1}=a$ for every $a\in A$,
$faf^{-1}\neq a$ for some $a\in A$ and $fAf^{-1}=A$. Let $Y_n$ be a
$n$-orbit of $A$ and $\overrightarrow{Y_n}$ be the arbitrary ordered set
$Y_n$, then $g\overrightarrow{Y_n}=\overrightarrow{Y_n}g$, $fY_n=Y_nf$, but
$f\overrightarrow{Y_n}\neq\overrightarrow{Y_n}f$.

\end{proposition}

\begin{lemma}\label{L_n}
Let $X_n$ be a $n$-orbit of a group $G$ and $L_n$ be a partition of $X_n$.
If the left action of $G$ on $X_n$ maintains $L_n$, then classes of $L_n$
are $n$-orbits of left cosets of some subgroup of $G$.

\end{lemma}
\proof
Let $Y_n\in L_n$. Since $L_n$ is a partition, the left action of
$Aut(Y_n)$ on $Y_n$ is transitive  and hence $Y_n$ is a $n$-orbit.
$\Box$\bigskip

The same we have

\begin{lemma}\label{R_n}
Let $X_n$ be a $n$-orbit of a group $G$ and $R_n$ be a partition of $X_n$.
If the right action of $G$ on $X_n$ maintains $R_n$, then classes of $R_n$
are $n$-orbits of right cosets of some subgroup of~$G$.

\end{lemma}

\paragraph{Intersections and unions of left- and right-automorphic
partitions.}

\begin{proposition}\label{Y_n*Z_n}
Let $Y_n$ and $Z_n$ be $n$-orbits, then $T_n=Y_n\cap Z_n$ is a $n$-orbit
and $Aut(T_n)=Aut(Y_n)\cap Aut(Z_n)$.

\end{proposition}
\proof
It is sufficient to choose an initial $n$-tuple from $T_n$.
$\Box$

\begin{definition}\label{covering}
Let $M$ be a set and $A$ be a system of subsets of $M$. We say that $A$ is
a covering of $M$ if classes of $A$ contain all elements of $M$ and have
non-vacuous intersections. If the all intersections are vacuous then we
say that (the covering) $A$ is a partition of $M$. So we say that $A$ is a
covering, if it is not a partition. We say also that $A$ is a covering, if
we do not know, whether it is a partition.

\end{definition}

\begin{definition}\label{automorphic}
Let $X_n$ be a $n$-orbit of a group $G$ and $Y_n$ be an arbitrary subset
of $X_n$, then we say that $L_n=GY_n$ is \emph{left-automorphic} and
$R_n=Y_nG$ is \emph{right-automorphic} covering of $X_n$.  This
definition we shall apply also to corresponding coverings of a $k$-orbit
$X_k$ of a group $G$ for $k<n$.

\end{definition}

\begin{corollary}\label{autom.-part.}
Let $L_n$ and $R_n$ be partitions  of a $n$-orbit of a group $G$ on
$n$-orbits of left and right cosets of a subgroup $A<G$, then $L_n$
and $R_n$ are left-automorphic and right-automorphic partitions.
\end{corollary}

\begin{definition}\label{part. op.}
Let $M$ be a set and $P,Q$ be partitions of $M$. We write:

\begin{itemize}
\item
$P\sqsubset Q$ if for every $A\in P$ there exists $B\in Q$, so that
$A\subset B$.

\item
$P\sqcap Q$ for partition of $M$ that consists of intersections of
classes from $P$ and $Q$.

\item
$P\sqcup Q$ for partition of $M$ whose class is a union of intersected
classes from $P$ and $Q$.

\end{itemize}
\end{definition}

\begin{proposition}\label{GA,GB;AG,BG}
Let $A,B<G$, then $GA\sqcap GB=G(A\cap B)$, $AG\sqcap BG=(A\cap B)G$,
$GA\sqcup GB=Ggr(A,B)$ and $AG\sqcup BG=gr(A,B)G$.

\end{proposition}
\proof
\begin{itemize}
\item
Since $GA\sqcap GB=G(GA\sqcap GB)$ and $A\cap B\in GA\sqcap GB$,
$GA\sqcap GB=G(A\cap B)$. Analogously $AG\sqcap BG=(A\cap B)G$.

\item
Let $\{A_i:\,i\in [1,l]\}\subset GA$, $\{B_j:\,j\in [1,m]\}\subset GB$,
$A_1=A$, $B_1=B$ and $U=\cup_{(i=1,l)} A_i=\cup_{(j=1,m)} B_j$, then
$U\in GA\sqcup GB$ and $e\in U$. Let $f\in A_i$, $g\in A_k$ and $B_j$ has
non-vacuous intersections with $A_i$ and $A_k$, then $f=gaba'$ for some
elements $a,a'\in A$ and $b\in B$. It shows that every element from $U$
can be represented as a product of elements from $A$ and $B$. Hence
$U=gr(A,B)$ and $GA\sqcup GB=Ggr(A,B)$.  The same $AG\sqcup
BG=gr(A,B)G$.~$\Box$

\end{itemize}

From this proposition follows:

\begin{lemma}\label{L_n,L_n;R_n,R_n}
Let $A,B<G$ and $X_n\in Orb_n(G)$, then $X_nA\sqcap X_nB=X_n(A\cap B)$,
$AX_n\sqcap BX_n=(A\cap B)X_n$, $X_nA\sqcup X_nB=X_ngr(A,B)$ and
$AX_n\sqcup BX_n=gr(A,B)X_n$.

\end{lemma}

\paragraph{Intersection and union of left-automorphic partition with
right-automorphic partition.}

Let $L_n=X_nA$ and $R_n=AX_n$. First we see that $L_n$ and $R_n$ have at
least one common class the $n$-orbit of $A$ containing initial $n$-tuple.
Then from proposition \ref{Yn=gYng-1->N} we know that if $L_n$ and $R_n$
have more as one common class, then $A$ has non-trivial normalizer in $G$.

\begin{lemma}\label{L_n,R_n,B<A}
Let $L_n\sqcap R_n$ be not trivial, i.e. it contains a class $Z_n$ by
power $l$, where $1<l<|A|$, then conjugate to $A$ subgroups have
non-trivial intersections and $Z_n$ is a $n$-orbit of a some subgroup
$B<A$.

\end{lemma}
\proof
Let $U_n\in L_n$, $W_n\in R_n$ and $Z_n=U_n\cap W_n$, then $U_n$ is a
$n$-orbit of some conjugate to $A$ subgroup $D$ and $W_n$ is a $n$-orbit
of $A$. Taking in opinion that we can choose an initial $n$-tuple from
$Z_n$, we obtain that $Z_n$ is a $n$-orbit of a subgroup $B=A\cap D$.
$\Box$

\begin{corollary}\label{L_n-sqcap-R_n,p}
Let $A$ be a prime order cyclic group, then $L_n\sqcap R_n$ is trivial.

\end{corollary}

The union $L_n\sqcup R_n$ can contain non-automorphic classes:
$$
\begin{array}{|c|}
\hline
123 \\
213 \\
\hline
132 \\
312 \\
\hline
231 \\
321 \\
\hline
\end{array}\ \sqcup\
\begin{array}{|c|}
\hline
123 \\
213 \\
\hline
132 \\
231 \\
\hline
312 \\
321 \\
\hline
\end{array}\ =\
\begin{array}{|c|}
\hline
123 \\
213 \\
\hline
132 \\
231 \\
312 \\
321 \\
\hline
\end{array}
$$
and therefore it is not of interest for investigation, nevertheless the
symmetry properties of this union can give an information about the
structure of the studied group $G$ and help to find subgroups of $G$ that
are supergroups for $A$.


\subsection{The actions of permutations on $k$-sets}

In order to consider the actions of permutations on $k$-sets we shall need
to have some special operations that we introduce from the beginning.


\subsubsection{Operations on $k$-sets}

\paragraph{Projecting and multiprojecting operators.}

Let $\alpha _k=\langle v_1\ldots v_k\rangle$ be a $k$-tuple, $l\leq
k$ and $i_1,i_2,\ldots ,i_l$ be $l$ different coordinates from $[1,k]$.
Then $\beta_l=\langle v_{i_1}\ldots v_{i_l}\rangle$ is a $l$-tuple that we
call a \emph{projection} of the $k$-tuple $\alpha_k$ on the ordered set of
coordinates $I_l=\{i_1<i_2<\ldots<i_l\}$. We shall enter a projecting
operator $\hat{p}$ and write this projection as
$\beta_l=\hat{p}(I_l)\alpha_k$. The projection of all $k$-tuples of a
$k$-set $X_k$ on $I_l$ gives a $l$-set $X_l=\hat{p}(I_l)X_k$.

The projection of all $k$-tuples of the $k$-set $X_k$ on $I_l$, that
distinguishes the equal $l$-tuples, is a multiset that we call a
\emph{multiprojection} of $X_k$ on $I_l$ and denote it as
$\uplus_{X_k}X_l$ or simply $\uplus X_l$, if from context it is clear,
what a multiprojection we consider. By definition $|\uplus X_l|=|X_k|$.
Using a multiprojecting operator $\hat{p}_{\uplus}$, we shall write that
$\uplus X_l=\hat{p}_{\uplus}(I_l)X_k$.

\paragraph{Concatenating operation.}

Let $\beta_l=\langle v_1\ldots v_l\rangle$ and $\gamma_m=\langle
u_1\ldots u_m\rangle$ be $l$- and $m$-tuple, then $(l+m)$-tuple
$\langle v_1v_2\ldots v_lu_1u_2\ldots u_m\rangle$ we call a concatenation
of $\beta_l$ and $\gamma_m$ and write this as $\beta _l\circ \gamma _m$.

It will be also suitable to use the concatenation of intersected tuples.
We shall consider such concatenation as multiset of coordinates, for
example $\langle 12\rangle\circ \langle 23\rangle=\langle 1223\rangle$.

We shall use this concatenating operation also for multisets of tuples in
the next way:

Let $\uplus Y_l$ and $\uplus Z_m$ be multisets with the same number of
tuples and $\phi:\uplus Y_l\leftrightarrow \uplus Z_m$, then $\uplus
Y_l\stackrel{\phi}{\circ}\uplus Z_m$ is the $(m+l)$-multiset, that
consists of concatenations of $l$-tuples of $Y_l$ with $m$-tuples of
$Z_m$ accordingly to the map $\phi$. We shall not write the map $\phi$, if
it is clear from context.

\paragraph{Operation properties.}\label{Op}

\begin{lemma}\label{pX_k}
$l$-projection of a $k$-orbit is a $l$-orbit.

\end{lemma}

Let $g$ be a permutation, $\alpha_n=\langle v_1\ldots v_n\rangle$ be a
$n$-tuple and $I_k$ be a $k$-subspace.

\begin{lemma}\label{pgan=gpan}
$g\hat{p}(I_k)\alpha_n=\hat{p}(I_k)g\alpha_n$.

\end{lemma}
\proof
It is sufficient to show the equality for $I_k\in [1,k]$.
$$g\hat{p}(I_k)\alpha_n=g\langle v_1\ldots v_k\rangle
                      =\langle gv_1\ldots gv_k\rangle$$ and
$$\hat{p}(I_k)g\alpha_n=\hat{p}(I_k)\langle gv_1\ldots gv_n\rangle
                      =\langle gv_1\ldots gv_k\rangle\!.\;\Box$$

\begin{lemma}\label{p.ang=pan.g}
$\hat{p}(I_k)(\alpha_ng)=\hat{p}(gI_k)\alpha_n$.

\end{lemma}
\proof
$$\hat{p}(I_k)(\alpha_ng)=\hat{p}([1,k])\langle v_{g1}\ldots
v_{gn}\rangle=\langle v_{g1}\ldots v_{gk}\rangle$$ and
$$\hat{p}(gI_k)\alpha_n=\hat{p}(g[1,k])\langle v_1\ldots v_n\rangle
                      =\hat{p}(\{g1<\ldots<gk\})\langle v_1\ldots
v_n\rangle =\langle v_{g1}\ldots v_{gk}\rangle.\ \Box $$

The equality
$\hat{p}(I_k)(\alpha_ng)=(\hat{p}(I_k)\alpha_n)g=\alpha_k(I_k,\alpha_n)g$
has not an interest application, because the corresponding right
permutation action on $k$-tuple $\alpha_k$ cannot be disengaged from the
$n$-tuple $\alpha_n$, as it takes place by the left permutation action on
$k$-tuple $\alpha_k$. But we shall write for convenience $\alpha_kg$
instead of $\hat{p}(I_k)(\alpha_ng)$, where it will not lead to
misunderstanding. Similarly, we consider a $k$-projection of equalities
$GY_n=X_nA$ and $AX_n=Y_nG$ (lemma \ref{GY_n=X_nA}).

From $\hat{p}(I_k)GY_n=\hat{p}(I_k)X_nA$, it follows
$G\hat{p}(I_k)Y_n=\hat{p}(AI_k)X_n$, where on definition \ref{GA,AG}
$$
  \hat{p}(AI_k)X_n\equiv \{\hat{p}(AI_k)\alpha_n:\,\alpha_n\in X_n\}
  \mbox{ and }\hat{p}(AI_k)\alpha_n=\{\hat{p}(aI_k)\alpha_n:\,a\in A\},
$$
so we can write the $k$-projection of this equality as
$GY_k(I_k)=\hat{p}(AI_k)X_n$ or (by correct understanding) simply as
$GY_k=X_kA$.

For the second equality we have
$\hat{p}(I_k)AX_n=A\hat{p}(I_k)X_n=AX_k(I_k)$ and
$$\hat{p}(I_k)Y_nG=\hat{p}(GI_k)Y_n\equiv \{\hat{p}(gI_k)Y_n:\, g\in G\}.$$
For convenience we will write the $k$-projection of second  equality simply
as $AX_k=Y_kG$.

\begin{proposition}\label{p(I_k)P_k*Q_k}
Let $X_n$ be a $n$-set and $P_n$, $Q_n$ be two partitions of $X_n$. Let
$P_k=\hat{p}(I_k)P_n$ and $Q_k=\hat{p}(I_k)Q_n$ be partitions of
$X_k=\hat{p}(I_k)X_n$. It does not necessitate the equality
$\hat{p}(I_k)(P_n\sqcap Q_n)=\hat{p}(I_k)P_n\sqcap \hat{p}(I_k)Q_n)$.

\end{proposition}
\proof
$X_n=\{15,26,36,45\}$, $P_n=\{\{15,26\},\{36,45\}\}$ and
$Q_n=\{\{15,36\},\{26,45\}\}$.~$\Box$


\subsubsection{Some additional definitions and auxiliary statements}

\begin{definition}
The presentation $GY_k=\{gY_k:g\in G\}$ that we use for the action of a
group $G$ on a $k$-subset $Y_k$ we shall apply also for the action of a
group $G$ on a system of $k$-subsets $P_k$ as $GP_k=\{gP_k:g\in G\}$ and
$gP_k=\{gY_k:Y_k\in P_k\}$.

\end{definition}
\begin{definition}
Let $M$ be a set and $Q$ be a set of subsets of $M$, then we write $\cup
Q\equiv \cup_{(U\in Q)}U$. From such definition follows that for $GP_k$ we
can consider two kinds of unions: $\cup GP_k$ and $\cup\cup GP_k$. For
example, let we have two sets $\{\{1,2\},\{2,3\}\}$ and
$\{\{2,3\},\{4,5\}\}$, then the first union of these sets is the set
$\{\{1,2\},\{2,3\},\{4,5\}\}$ and the second is $\{1,2,3,4,5\}$.

The symbol $\sqcup$ we shall apply to union of intersected classes of a
system of sets, as in example: $\sqcup
\{\{1,2\},\{2,3\},\{4,5\}\}=\{\{1,2,3\},\{4,5\}\}$.

\end{definition}

\begin{definition}
Let $\alpha_k$ be a $k$-tuple. We shall write the set of coordinates of
$\alpha_k$ as $Co(\alpha_k)$. We shall use this notation also for
a $k$-set $Y_k$, where $Co(Y_k)=\{Co(\alpha_k):\alpha_k\in Y_k\}$.

\end{definition}

\begin{lemma}\label{rm}
Let $M$ be a $m$-element set and $\uplus M$ be a homogeneous multiset
with multiplicity $k$. Let $P_m$ be a multipartition of $\uplus M$ on
$k$-element multisubsets of $M$, then $P_m$ is a union of $k$
distributions of $m$ elements of set $M$ between $m$ $k$-element classes
of $P_m$.

\end{lemma}
\proof
We do an induction on $m$. Let us to represent $P_m$ as $m\times k$
matrix, whose lines are classes of $P_m$. Let $m=1$, then the statement is
evidently correct. Let $m>1$ and first $l\leq k$ lines contains all $k$
occurences of an element $u\in M$. By permutation of elements in these $l$
lines we can placed element $u$ in all $k$ columns. Now by permutation of
elements in $k$ columns we can replaced element $u$ in the first line.
Thus we have obtained $(m-1)\times k$ matrix (without the first line) that
by induction hypothesis can be transformed (with permutation of elements
in lines) to $m-1$ different elements in each column. Now we need only to
do the inverse permutation of element $u$ from the first line with
corresponding elements in other $l-1$ lines. $\Box$\bigskip

From this lemma follows

\begin{corollary}\label{mMmM->MM}
Let $M$ be a $m$-element set and $M_2=\uplus M\circ \uplus M$, where
$\uplus M$ is homogeneous, then $M_2$ can be partition in $m$-element
subsets that are concatenations $M\circ M$.

\end{corollary}
\proof
Let $|\uplus M|/|M|=k$, $M_2$ be associated with space $I_2=I_1^1\circ
I_1^2$ and $L_2$ be the partition  of $M_2$ on $k$-element classes so that
$\hat{p}(I_1^1)L_2=M$, then $\hat{p}(I_1^2)L_2=P_m$ from lemma
\ref{rm}. $\Box$

\begin{lemma}\label{mMmM}
Let $M_2$ be a set of pairs that is associated with space
$I_2=I_1^1\circ I_1^2$. Let $\hat{p}_{\uplus}(I_1^1)M_2=
\hat{p}_{\uplus}(I_1^2)M_2$, then $M_2$ can be partition in cycles.

\end{lemma}
\proof
Let $\langle u_1u_2\rangle\in M_2$, then there exists a pair $\langle
u_2u_3\rangle\in M_2$.  The continuation gives a first cycle $C_2=\{\langle
u_1u_2\rangle,\langle u_2u_3\rangle,\ldots, \langle u_ru_1\rangle\}$. The
set $M_2\setminus C_2$ holds the property of the set $M_2$.
$\Box$\bigskip

\begin{definition}
Let $Y_k$ be a $k$-subset defined on the subset $U\subset V$, then under
$Aut(Y_k)$ we shall understand the permutation group on the set $U$. An
extension of $Aut(Y_k)$ on the set $V$ we shall write as $Aut(Y_k;V)$.

\end{definition}

\begin{definition}\label{sdeco}
Let $c_1c_2\ldots c_l$ be decomposition of a permutation $g$. Product $a$
of some cycles from this decomposition we shall call a {\em
subdecomposition} of $g$ and write this as $a\subset g$.

Let $a$ be a subdecomposition of an automorphism $g\in G$. We shall call
the permutation $a$ as a \emph{subautomorphism} of $G$ and write this as
$a\subset\in G$. Let $B$ be an intransitive subgroup of $G$ and $A(U)$ be
a transitive component of $B$ on the subset $U\subset V$. We shall write
this fact as $A(U)<<G$ and say that $A(U)$ is a projection of $B$. It is
clear that $A(U)$ is generated by some subautomorphisms of $G$.

We can consider an action of a subautomorphism $a$ on a $k$-orbit $X_k$ of
a group $G$, extending it to an action of some automorphism $g\in G$.

\end{definition}

\begin{definition}
Let $X_k$ be a $k$-orbit of a group $G$ and $Y_k\subset X_k$ be a
$k$-orbit of a subgroup of $G$, then we say that $Y_k$ is a $k$-suborbit
of $G$.

\end{definition}

\begin{definition}
Let $G$ be a group, $X_n\in Orb_n(G)$, $I_k$ be a $k$-subspace and
$Co(I_k)$ be a $1$-orbit of some subgroup $A<G$, then we say that number
$k$ is automorphic, $I_k$ is an \emph{automorphic subspace} and
$X_k=\hat{p}(I_k)X_n$ is a \emph{right-automorphic $k$-orbit} or a
\emph{$k$-rorbit}.

\end{definition}


\subsubsection{The left action of permutations on $k$-sets}

The left action of a permutation on a $k$-set is the same as its action
on a $n$-set. The right action of a permutation on a $k$-set follows
from its action on a $n$-set. The right action is not just visible
combinatorially as the left action, so we shall begin with the left
action.

We say that two sets of $k$-tuples $Y_k$ and $Y_k'$ are {\em
$S_n$-isomorphic} or simply isomorphic if there exists permutation $g\in
S_n$ so that $Y_k'=gY_k$, for example $Y_2=\{\langle 12\rangle,\langle
21\rangle\}$ and $Y_2'=\{\langle 13\rangle,\langle 31\rangle\}$. We shall
say that $Y_k$ and $Y_k'$ are $G$-isomorphic, if $g\in G<S_n$ and we study
invariants of $G$.  We shall not indicate a group relative to that we
consider the symmetry, if it is clear from context. From this definition
follows:

\begin{proposition}
Let  $X_k$ be a $k$-orbit and $Y_k$ be an arbitrary $k$-subset of $X_k$,
then $Aut(X_k)Y_k$ is a covering of $X_k$ on isomorphic to $Y_k$ classes.

\end{proposition}

and

\begin{corollary}
$n$-Orbits of left cosets of a subgroup $A$ of a group $G$ are isomorphic.

\end{corollary}

The $n$-orbits of right cosets of a subgroup $A$ in general are not
isomorphic. An example is $3$-orbits of the subgroup $A=gr((12))<S_3$:

$$
\begin{array}{||ccc||}
\hline
1 & 2 & 3 \\
2 & 1 & 3 \\
\hline
\end{array}\
\mbox{ and }\
\begin{array}{||ccc||}
\hline
1 & 3 & 2 \\
2 & 3 & 1 \\
\hline
\end{array}\,.
$$
The same is valid for $k$-orbits of left, right cosets of a subgroup $A$ by
$k<n$.

\begin{proposition}
Let $G$ be a group, $X_n\in Orb_n(G)$, $X_k\in Orb_k(G)$ and $A<G$, then
\begin{itemize}
\item
$n$-Orbits of left cosets of $A$ form a partition of $X_n$.
\item
The $G$-isomorphic $k$-orbits of left cosets of $A$ belong to the same
$k$-orbit $X_k$ and form a covering of $X_k$.

\end{itemize}
\end{proposition}
\proof
The first statement is evident. Let $Y_k\in Orb_k(A)$ be subset of $X_k$ ,
then a covering $L_k=GY_k$ of $X_k$ contains all $G$-isomorphic to $Y_k$
$k$-orbits. The example of such covering is $L_1=\{\{1,2\}\{2,3\}\{1,3\}\}$
$\Box$\bigskip

A $k$-orbit $X_k$ of a group $G$ can have different representations through
$k$-orbits of the same subgroup $A<G$, because $X_k$ can contain
non-isomorphic $k$-orbits of $A$. For example, $1$-orbit of the symmetric
group $S_n$ can be represented, on the one hand, as a covering by
$1$-orbits of left cosets of $gr((12\ldots (n-1)))$ that are
$(n-1)$-element subsets of $V$ and, on the other hand, as a partition on
$1$-orbits of left cosets of this subgroup that are $1$-element subsets of
$V$.

\begin{lemma}\label{Yk->Un}
Let $A<G$,  $X_n\in Orb_n(G)$, $Y_n\subset X_n$ be a $n$-orbit of $A$ and
$Y_k=\hat{p}(I_k)Y_n$.  Let $U_n\subset X_n$ be the union of such classes
of $GY_n$ whose $I_k$-projection is $Y_k$, then $U_n$ is a $n$-orbit.

\end{lemma}
\proof
$Aut(U_n)=Aut(Y_k;V)\cap G$. $\Box$\bigskip

The subset $U_n$ can contain not all $n$-tuples whose $I_k$-projections
belongs to $Y_k$. An example is given by the group $S_3$,
$U_3=Y_3=\{123,213\}$ and $k=1$:

$$
\begin{array}{||c|cc||}
\hline
1 & 2 & 3 \\
2 & 1 & 3 \\
\hline
1 & 3 & 2 \\
3 & 1 & 2 \\
\hline
2 & 3 & 1 \\
3 & 2 & 1 \\
\hline
\end{array}\,.
$$
In this case $U_1=\{1,2\}$ and the intersections $\{1,2\}\cap \{1,3\}$ and
$\{1,2\}\cap \{2,3\}$ are not vacuous.

\paragraph{Now we consider when a subgroup $A<G$ forms a partition of a
$k$-orbit $X_k\in Orb_k(G)$.}

Let a $k$-set $X_k$ be a $k$-projection of a $(k+1)$-set $X_{k+1}$, then we
shall say that $X_{k+1}$ is an \emph{extension} of $X_k$ on a
$(k+1)$-subspace or is a $(k+1)$-extension of $X_k$.

\begin{lemma}\label{Xk-X(k+1)}
Let $X_k\in Orb_k(G)$, $A<G$ and $X_kA=GY_k$ be a partition. Let
$X_{k+1}\in Orb_{k+1}(G)$ be an extension of $X_k$, then $X_kA$ generates
a partition $X_{k+1}A$.

\end{lemma}
\proof
Let $Y_{k+1}\subset X_{k+1}$ and $Y_k=\hat{p}(I_k)Y_{k+1}$. If $GY_{k+1}$
is not a partition, then evidently $GY_k$ contains intersected classes.
$\Box$

\begin{lemma}\label{XAUXB=XAB}
Let $X_k\in Orb_k(G)$, $A,B<G$ be not conjugate subgroups or $A,B$ be
conjugate subgroups and $Y_k,Z_k\subset X_k$ be not isomorphic
$k$-orbits of $A$ and $B$ correspondingly. Let $Y_k,Z_k$ be
intersected and $X_kA=GY_k,X_kB=GZ_k$ be partitions, then
$X_kA\sqcup X_kB=X_kgr(A,B)$.

\end{lemma}
\proof
$G(GY_k\sqcup GZ_k)=GY_k\sqcup GZ_k$. Then $X_kA\sqcup X_kB=GT_k=X_kC$,
where $C=gr(A,B)$ accordingly to lemma \ref{L_n,L_n;R_n,R_n}.
$\Box$\bigskip

Without condition $Y_k\cap Z_k\neq\emptyset$ the equality $X_kA\sqcup
X_kB=X_kgr(A,B)$ can lead to misunderstanding for conjugate subgroups
$A,B$. Consider an example:

$$
\begin{array}{||cc||}
\hline
1 & 2  \\
2 & 1  \\
\hline
1 & 3  \\
3 & 1  \\
\hline
2 & 3  \\
3 & 2  \\
\hline
\end{array}\,,\
\begin{array}{||cc||}
\hline
1 & 3  \\
2 & 3  \\
\hline
1 & 2  \\
3 & 2  \\
\hline
2 & 1  \\
3 & 1  \\
\hline
\end{array}\,.
$$
\begin{itemize}
\item
The subgroups $A=gr((12))$ and $B=gr((13))$ are conjugate in $G=S_3$, so
they determine the same partitions of $2$-orbits of $G$, if we
consider isomorphic (not intersected) $2$-orbits of these subgroups.
Therefore in this case $X_kA\sqcup X_kB=X_kA\neq X_kgr(A,B)$.

\item
If we consider not isomorphic and not intersected $2$-orbits of the same
subgroup $A=gr((12))$, then we have two different partitions $L_1,L_2$ of
$2$-orbit $X_2$ of $S_3$. So it can be seen that $X_2A\sqcup
X_2A\neq X_2gr(A,A)=X_2A$. This misanderstansing follows from
interpretation $X_2A$, first, as $S_3Y_2'$, $Y_2'=\{12,21\}$, and, second,
as $S_3Y_2''$, $Y_2''=\{13,23\}$, where $Y_2'\cap Y_2''=\emptyset$. If we
take $Y_2''=\{12,32\}$, then $A=gr((12))$, $B=gr((13))$ and formula
$X_kA\sqcup X_kB=X_kgr(A,B)$ gives correct result.

\end{itemize}

\begin{lemma}\label{Co(Y_k)=Co(alpha_k)}
Let $X_k\in Orb_k(G)$, $Y_k\subset X_k$, $\alpha_k\in Y_k$ and $Y_k$ be a
maximal subset with property: $Co(Y_k)=Co(\alpha_k)$, then $L_k=GY_k$ is a
partition of $X_k$.

\end{lemma}
\proof
On the condition the permutations of coordinates of $Y_k$ that
maintains $Y_k$ maintains also each class of $L_k=GY_k$. $\Box$

\paragraph{In the next statements we shall study reverse question.}
Namely, when a subset of a $k$-orbit of a group $G$ generates a subgroup
of $G$.

\begin{lemma}\label{L_k}
Let $X_k\in Orb_k(G)$ and $L_k$ be a partition of $X_k$. If the left
action of $G$ on $X_k$ maintains $L_k$, then classes of $L_k$ are
$k$-orbits of left cosets of some subgroup $A<G$.

\end{lemma}
\proof
Let $Y_k\in L_k$, then $L_k=GY_k$, hence a subgroup $A=Aut(Y_k:V)\cap G$
acts on $Y_k$ transitive. $\Box$

\begin{remark}
It can be seen that the partitioning of $X_n$ (and hence $X_k$) on
$G$-isomorphic classes is not sufficient for the automorphism of this
partition.  This shows the next partition $L_n$ of a $n$-orbit $X_n$ of
the group $G=C_6$:

$$
\begin{array}{||cccccc||}
\hline
1 & 2 & 3 & 4 & 5 & 6 \\
2 & 3 & 4 & 5 & 6 & 1 \\
\hline
3 & 4 & 5 & 6 & 1 & 2 \\
4 & 5 & 6 & 1 & 2 & 3 \\
\hline
5 & 6 & 1 & 2 & 3 & 4 \\
6 & 1 & 2 & 3 & 4 & 5 \\
\hline
\end{array}\, .
$$
The classes of $L_n$ are connected with permutation $(135)(246)$,
but are not automorphic.  Let $Y_n\in L_n$, then in this case $|GY_n|$ is
a covering of $X_n$.

\end{remark}

From lemmas \ref{Co(Y_k)=Co(alpha_k)} and \ref{L_k} follows

\begin{corollary}\label{Y_k-A}
Let $X_k\in Orb_k(G)$, $Y_k\subset X_k$, $\alpha_k\in Y_k$ and $Y_k$ be a
maximal subset with property: $Co(Y_k)=Co(\alpha_k)$, then $Y_k$ is a
$k$-suborbit of $G$.

\end{corollary}

Consider a generalization of lemma \ref{L_k}.

\begin{theorem}\label{|Y_k||GY_k|}
Let $X_k\in Orb_k(G)$, $Y_k\subset X_k$ and $L_k=GY_k$ be a covering of
$X_k$.  If $|Y_k||L_k|$ divides $|G|$, then $Y_k$ is a $k$-orbit of some
subgroup $A<G$.

\end{theorem}
\proof
On condition, there exists a partition $L_n$ of a $n$-orbit $X_n\in
Orb_n(G)$ so that $\hat{p}(I_k)L_n=L_k$ and $|L_n|=|L_k|$. Then there
exists $Y_n\in L_n$ so that $\hat{p}(I_k)Y_n=Y_k$,
$G\hat{p}(I_k)Y_n=\hat{p}(I_k)GY_n=GY_k$. Hence $L_n=GY_n$ is a partition
of $X_n$ and we can apply the lemma \ref{L_n}. $\Box$

\begin{corollary}\label{Qk}
Let $X_k\in Orb_k(G)$, $Y_k\subset X_k$ and $L_k=GY_k$ be a covering of
$X_k$. Let classes of $L_k$ can be assembled in isomorphic partitions of
$X_k$.  Let $Q_k$ be a system of these partitions and $|Q_k||X_k|$ divides
$|X_n|$. Let $L_k'\in Q_k$ and $Y_k\in L_k'$, then $X_k$ is a $k$-orbit of
a subgroup $B=Aut(L_k')\cap G$, $L_k'=BY_k$ and $Y_k$ is a $k$-orbit of
some subgroup $A<B$.

\end{corollary}
\proof
Let $X_n\in Orb_n(G)$ and $X_k=\hat{p}(I_k)X_n$. Since $|Q_k||X_k|$
divides $|X_n|$, there exists a partition $L_n$ of $X_n$ so that
$\hat{p}(I_k)L_n=X_k$ and permutations of classes of $L_n$ correspond with
permutations of partitions from $Q_k$. So the classes of $L_n$ are
$n$-orbits of subgroups of $G$ and hence $X_k$ is a $k$-orbit of the
subgroup $B=Aut(L_k')\cap G$. Since $BY_k=L_k'$ is a partition of $X_k$
and $B$ acts on $L_k'$ transitive, the subgroup $A=Aut(Y_k;V)\cap B$ acts
on $Y_k$ transitive. $\Box$\bigskip

The example of such $G$-isomorphic system of partitions is

$$
Q_1=
\{\
\begin{array}{||c||}
\hline
1  \\
2  \\
\hline
3  \\
4  \\
\hline
\end{array}\,,\
\begin{array}{||c||}
\hline
1  \\
3  \\
\hline
2  \\
4  \\
\hline
\end{array}\,,\
\begin{array}{||c||}
\hline
1  \\
4  \\
\hline
2  \\
3  \\
\hline
\end{array}\
\},
$$
where  $Q_1$ is formed by $1$-orbits of left cosets of the subgroup
$gr((12)(34))<A_4$.

Let $Q_k$ be a set of $S_n$-isomorphic partitions that therefore do not
belong to the same $k$-orbit $X_k$, then the action of $G$ on $X_k$
maintains simultaneously all partitions $L_k^i$ as in the previous
example, where now $Q_1$ is formed by $1$-orbits of not conjugate
subgroups $gr((12)(34))$, $gr((13)(24))$ and $gr((14)(23))$ of the group
$S_2\otimes S_2$.

\begin{corollary}\label{rCk-gr}
Let $X_k$ be a $k$-rorbit of a group $G$ that contains a $k$-rcycle
$rC_k$, then $rC_k$ is a $k$-orbit of some subgroup $A<G$.

\end{corollary}
\proof
It is a special case of the corollary \ref{Y_k-A}. $\Box$\bigskip

The order of the subgroup $A$ in the corollary can differ from $k$. The
example gives the subgroup $A=gr((1234)(56))<S_6$ and $rC_2=\{\langle
56\rangle,\langle 65\rangle\}$.

Projections of $k$-rcycles from $X_k$ on $l$-subspace ($l<k$) can have
non-trivial intersections as in example:

$$
\begin{array}{lcr}
\begin{array}{||cc|c||}
\hline
1 & 2 & 3  \\
\hline
2 & 3 & 1  \\
3 & 1 & 2  \\
\hline
\end{array}&
\mbox{ and }&
\begin{array}{||cc|c||}
\hline
1 & 2 & 4  \\
\hline
2 & 4 & 1  \\
4 & 1 & 2  \\
\hline
\end{array}\,.
\end{array}
$$
So these projections form a covering of $X_l$.


\subsubsection{The special left action of permutations on $k$-sets}

There exists a left action of permutations on a $k$-orbit $X_k$ of a group
$G$ that forms a partition $R_k$ of $X_k$ on $k$-orbits of right cosets of
a subgroup $A<G$. It is a partition $R_k=AX_k$.

Classes of $AX_k$ as well as classes of $AX_n$ are in general case not
isomorphic. Moreover, if the classes of $AX_n$ have the same order, then
the classes of $AX_k$ satisfy to this property not always. For example
$R_1=gr((12))\{1,2,3\}=\{\{1,2\},\{3\}\}$.

$k$-Orbits of left cosets of a subgroup $A$  can have intersections,
$k$-Orbits of right cosets of a subgroup $A$ have no intersection.

$k$-Orbits of left cosets of a subgroup $A$ are $k$-orbits of subgroups
that are conjugate to $A$, $k$-orbits of right cosets of a subgroup $A$
are $k$-orbits of the subgroup $A$. These properties we shall assemble in
the following statements:

\begin{lemma}\label{conjugate-n}
Let $A,B$ be conjugate subgroups of a group $G$ and $X_n$ be a $n$-orbit
of $G$, then

\begin{enumerate}
\item
The partitions of $G$ on left (right) cosets of subgroups $A,B$ are not
equal.

\item
The partitions of $X_n$ on $n$-orbits of left cosets of subgroups $A,B$
are equal and this partition consists of isomorphic classes.

\item
The partitions of $X_n$ on $n$-orbits of right cosets of subgroups $A,B$
are not equal and each partition consists of not isomorphic classes of
power $|A|$.

\end{enumerate}
\end{lemma}
\proof
The first statement is the fact from the group theory, the second is the
repeating of lemma \ref{gY_n} and the third follows from lemma
\ref{Y_ng}.  $\Box$

\begin{corollary}\label{conjugate-k}
Let $A,B$ be conjugate subgroups of a group $G$ and $X_k$ be a $k$-orbit
of $G$, then

\begin{enumerate}
\item
The coverings of $X_k$ on isomorphic $k$-orbits of left cosets of
subgroups $A,B$ are equal.

\item
The partitions of $X_k$ on $k$-orbits of right cosets of subgroups $A,B$
are not equal, each partition consists of not isomorphic classes, which can
differ by power.

\end{enumerate}
\end{corollary}

\paragraph{$k$-orbit property of normal subgroups.}

\begin{lemma}\label{AX_k=X_kA}
Let $X_k\in Orb_k(G)$, $A<G$ and $R_k=AX_k=X_kA=L_k$. If $A$ is the
maximal subgroup, then $A\lhd G$.

\end{lemma}
\proof
Let $X_n\in Orb_n(G)$, then $R_n=AX_n=X_nA=L_n$, because $L_n$ is the only
partition of $X_n$ with property $\hat{p}(I_k)L_n=L_k$ and $|L_n|=|L_k|$,
hence $A\lhd G$. $\Box$\bigskip

\begin{corollary}\label{fix-tuple}
Let $Y_k$ be a maximal subset of a $k$-orbit $X_k$ so that
$Co(Y_k)=Co(\alpha_k)$ for some $k$-tuple $\alpha_k\in Y_k$, then a
stabilizer $Stab(\alpha_k)\lhd Stab(Y_k)$.

\end{corollary}

\begin{lemma}\label{A<|G}
Let $A\lhd G$ and $X_k\in Orb_k(G)$, then $R_k=AX_k=X_kA=L_k$.

\end{lemma}
\proof
It is given that $R_n=AX_n=X_nA=L_n$, then
$\hat{p}(I_k)AX_n=\hat{p}(I_k)X_nA$ or $AX_k=X_kA$.~$\Box$\bigskip

So we have

\begin{theorem}\label{s-gr}
A group $G$ is a simple group if and only if $AX_k\neq X_kA$ for
arbitrary $k$, arbitrary $X_k\in Orb_k(G)$ and each subgroup $A<G$.

\end{theorem}

\paragraph{Intersections and unions of $k$-orbits}

Above we have seen (proposition \ref{Y_n*Z_n}) that the intersection of
$n$-orbits is a $n$-orbit. The same is correct

\begin{lemma}\label{Y_k*Z_k}
Let $Y_k$ and $Z_k$ be $k$-orbits, then $T_k=Y_k\cap Z_k$ is a $k$-orbit
of $Aut(T_k)=Aut(Y_k)\cap Aut(Z_k)$.

\end{lemma}
\proof
It is sufficient to consider the intersection of $n$-orbits of $Aut(Y_k)$
and $Aut(Z_k)$ and then their corresponding $k$-projections.
$\Box$\bigskip

\begin{corollary}\label{R_k*R_k'}
Let $Y_k$ and $Z_k$ be $k$-suborbits of a $k$-orbit $X_k$, $A=Aut(Y_k)\cap
Aut(X_k)$ and $B=Aut(Z_k)\cap Aut(X_k)$, then $(A\cap B)X_k=AX_k\sqcap
BX_k$.

\end{corollary}

For subgroups $G<Aut(Y_k)$ and $G'<Aut(Z_k)$ with $k$-orbits $Y_k\in
Orb_k(G)$ and $Z_k\in Orb_k(G')$ the similar relation is not correct. Let
us to give a counterexample:

The partition $R_1=AX_1=AV$ is a system of orbits of $A$ on $V$ in the
conventional meaning. Let $A=gr((12))<S_3$ and $B=gr((123))<S_3$, then
$AV=\{\{1,2\},\{3\}\}$, $BV=\{\{1,2,3\}\}$, $AV\sqcap
BV=\{\{1,2\},\{3\}\}$ and $(A\cap B)V=\{\{1\},\{2\},\{3\}\}$.

This example, lemma \ref{Y_k*Z_k} and corollary \ref{R_k*R_k'} determine
the relation between intersections of groups, their $k$-orbits and
corresponding systems of $k$-orbits of right cosets of this groups
(subgroups). The same is valid for the intersection of systems of
$k$-orbits of left cosets.

The union of partitions of $X_k\in Orb_k(G)$ on $k$-orbits of left cosets
of subgroups $A,B<G$ we have considered in lemma \ref{XAUXB=XAB}.

For the union of partitions of $X_k\in Orb_k(G)$ on $k$-orbits of right
cosets of subgroups $A,B<G$ we have:

\begin{lemma}\label{AV*BV,AV+BV}
$AX_k\sqcup BX_k=gr(A,B)X_k$.

\end{lemma}

\paragraph{The condition of automorphism of a subspace $I_k\subset V$.}

\begin{lemma}\label{k-rorbit}
Let $X_n$ be a $n$-orbit of a transitive group $G$, $X_k=\hat{p}(I_k)X_n$,
$|X_n|/|X_k|>1$ and $I_k$ contains all elements of $V$ that are fixed with
$Stab(I_k)<G$, then $I_k$ is automorphic.

\end{lemma}
\proof
Let $Co(I_k)=\{v_i:\, i\in [1,k]\}$, $T_k^1\in Orb_k(Stab(v_1))$,
$T_k^1\subset X_k$, $\hat{p}(v_1)T_k^1=v_1$ and $L_k'=GT_k^1$. Let
$T_k^i\in Orb_k(Stab(v_i))$ be a class of $L_k'$, then
$Stab(I_k)T_k^i=T_k^i$. It follows that all $T_k^i$ consist of $k$-orbits
of right cosets of $Stab(I_k)$ and systems of these $k$-orbits of right
cosets of $Stab(I_k)$ for different $i$ are isomorphic. Hence each $T_k^i$
contains a fix $k$-tuple $\alpha_k^i$ for that $Co(\alpha_k^i)=Co(I_k)$.
The union $\cup_{(i=1,k)} \alpha_k^i$ is a $k$-orbit of a normalizer
$N_G(Stab(I_k))$ (corollary \ref{fix-tuple}) that acts transitive on the
subset $Co(I_k)$. Hence $I_k$ is automorphic.  $\Box$\bigskip


\subsubsection{The right action of permutations on $k$-sets}
\paragraph{Right action isomorphism.}
Under right action of a permutation $g$ on a $k$-tuple $\alpha_k$ we
understand the mapping of $\alpha_k$ on a $k$-tuple $\beta_k\equiv
\alpha_kg$ that is placed on the position of coordinates of $\alpha_k$ in
the same $n$-tuple $\alpha_n$ under the right action of $g$ on $\alpha_n$.
If $\alpha_k=\langle v_1v_2\ldots v_k\rangle$, then on definition we write
$\alpha_kg=\langle v_{g1}v_{g2}\ldots v_{gk}\rangle$. Thus we consider the
$k$-tuple $\alpha_k\subset \alpha_n$ with its certain position in
$\alpha_n$ that we define by $k$-subspace of coordinates $I_k$.

Let $X_n\in Orb_n(G)$, $Y_n\subset X_n$, $Y_k=\hat{p}(I_k)Y_n$, $g\in S_n$
and $Y_k'=Y_kg=\hat{p}(gI_k)Y_n$, then we say that $Y_k$ and $Y_k'$ are
right $S_n$-related. If $g\in G$, then $Y_k$ and $Y_k'$ are right
$G$-related. In general case the image $Y_k'$ is not isomorphic to its
original $Y_k$. We shall study, when the right action of a permutation
transforms $k$-subset $Y_k$ on isomorphic $k$-subset $Y_k'$.

\paragraph{Two kinds of right action isomorphism.}
\begin{definition}\label{bYk=Yka}
If we study $k$-orbits of a group $G$, $X_k\in Orb_k(G)$, $Y_k\subset
X_k$, $a,b\in S_n$ and a $k$-subset $Y_k'=Y_ka=bY_k$, then we say that
$Y_k'$ and $Y_k$ are right $S_n$-isomorphic (as in first case of example
\ref{simple-struc}). If $a,b\in G$, then we say that $Y_k'$ and $Y_k$ are
right $G$-isomorphic.

\end{definition}

Let $Y_k=\hat{p}(I_k)Y_n\subset X_k$, $X_k'\in Orb_k(G)$,
$Y_k'=\hat{p}(I_k')Y_n\subset X_k'$ and $Y_k'=Y_ka=bY_k$. If $a\in G$ then
$X_k'=X_k$, and $Y_k'\subset X_k$ too. Now we shall find, when from $a\in
G$ it follows $b\in G$.

\begin{proposition}\label{SnYk-not-Aut(Xk)Yk}
Let $Y_k$ and $Y_k'$ be arbitrary $S_n$-isomorphic subsets of a $k$-orbit
$X_k$, then $Y_k$ and $Y_k'$ are not with necessary $Aut(X_k)$-isomorphic.

\end{proposition}
\proof
The $2$-subsets $\{\langle 12\rangle,\langle 23\rangle\}$ and $\{\langle
12\rangle,\langle 24\rangle\}$ from the $2$-orbit $X_2'$  (page
\pageref{X_2'}) are $S_n$-isomorphic, but not $Aut(X_2')$-isomorphic.
$\Box$\bigskip

\begin{proposition}\label{SnYk->Aut(Xk)Yk}
Let $Y_k$ and $Y_k'$ be arbitrary right $S_n$-isomorphic $k$-subsets of
a $k$-orbit $X_k$ of a group $G$, then $Y_k$ and $Y_k'$ are
$G$-isomorphic.

\end{proposition}
\proof
Let $Y_k=\{I_k(i):i\in [1,|Y_n|]\}$ and $Y_k'=\{I_k'(i):i\in [1,|Y_n|]\}$,
then $Y_k=\hat{p}(I_k(i))Y_n$ and $Y_k'=\hat{p}(I_k'(i))Y_n$. Maps
$I_k(i)\rightarrow I_k'(i)$ are restrictions of automorphisms. Hence a map
$\cup I_k(i)\rightarrow \cup I_k'(i)$ is also a restriction of an
automorphism. $\Box$\bigskip

\paragraph{General properties of right permutation action.}

\begin{proposition}\label{Yk'=gYk''f-1}
$k$-Orbits of left, right cosets of a subgroup $A$ are connected with
elements of~$G$.

\end{proposition}
\proof
From proposition \ref{Yn'=gYn''f-1} we obtain
$Y_k'(I_k)=\hat{p}(I_k)Y_n'= \hat{p}(I_k)gY_n''f^{-1}
=g\hat{p}(I_k)Y_n''f^{-1} = gY_k''(I_k)f^{-1} = gY_k''(f^{-1}I_k)$.
$\Box$\smallskip

\begin{proposition}\label{H->r.aut}
The right action of any automorphism $g\in G$ on a $k$-orbit of a normal
subgroup $H\lhd G$ is isomorphic.

\end{proposition}
\proof
Accordingly to proposition \ref{H->Yn=gYng-1}
$\hat{p}(I_k)gY_n=\hat{p}(I_k)Y_ng$ or $gY_k=Y_kg$. $\Box$\bigskip

Now we give a generalization of proposition \ref{Yn=gYng-1->N}.

\begin{lemma}\label{Lk*Rk}
Let $A<G$, $X_k\in Orb_k(G)$, $L_k=X_kA=GY_k$, $R_k=AX_k$ and $Q_k=L_k\cap
R_k\neq \emptyset$. Let $A$ be a maximal subgroup, then $A$ has non-trivial
normalizer $B=N_G(A)$.

\end{lemma}
\proof
Let $Z_k=\cup Q_k$, then $Aut(Y_k;V)\cap Aut(Z_k;V)\lhd Aut(Z_k;V)$, where
$Z_k$ is not with necessary automorphic. Hence $A=Aut(Y_k;V)\cap
Aut(Z_k;V)\cap G\lhd Aut(Z_k;V)\cap G=B$.  $\Box$

\begin{proposition}\label{gakg-1,fYkf-1}
Let $A$ be a subgroup of $G$, $g,f\in G$, $gag^{-1}=a$ for every $a\in A$
and $faf^{-1}\neq a$ for some $a\in A$, but $fAf^{-1}=A$. Let $Y_k$ be a
$k$-orbit of $A$ and $\overrightarrow{Y_k}$ be the arbitrary ordered set
$Y_k$, then $g\overrightarrow{Y_k}=\overrightarrow{Y_k}g$ and $fY_k=Y_kf$,
but $f\overrightarrow{Y_k}\neq\overrightarrow{Y_k}f$.

\end{proposition}
\proof
It follows from proposition \ref{g*an*g-1,f*Yn*f-1}. $\Box$

\paragraph{Right permutation action on $k$-rcycles.}
\begin{lemma}\label{C_k-conc}
Let $C_{lk}$ be a $(l,k)$-cycle, then $C_{lk}=\uplus C_{l_1p_1}\circ
\uplus C_{l_2p_2}\circ\ldots \circ \uplus C_{l_qp_q}$, where
$\sum_{i=1}^qp_i=k$, $l_i$ divides  $|\uplus C_{l_ip_i}|=l$, $C_{l_ip_i}$
is a $p_i$-projection of a $l_i$-rcycle and different $l_i$-rcycles have
no intersection on $V$.

\end{lemma}
\proof
From definition it follows that $C_{lk}=gr(g)\alpha_k$ for some
permutation $g$ and $k$-tuple $\alpha_k$. Let $n$-tuple
$\beta_n=\beta_{l_1}\circ \beta_{l_2}\circ \ldots\circ \beta_{l_q}$,
$g=(\beta_{l_1})(\beta_{l_2})\ldots (\beta_{l_q})$ and
$m=\sum_{i=1}^ql_i$, then $g$ generates $(l,m)$-cycle $C_{lm}=\uplus
rC_{l_1}\circ \uplus rC_{l_2}\circ\ldots\circ \uplus rC_{l_q}$. The
$(l,k)$-cycle $C_{lk}$ is a $k$-projection of the $(l,m)$-cycle
$C_{lm}$.~$\Box$\bigskip

The $(l,k)$-cycle $C_{lk}$ can be represented as a concatenation of
$(l,p_i)$-multiorbits, whose $p_i$-projections are either $p_i$-tuple, or
$S_{l_i}^{p_i}$-orbits, or $p_i$-projections of $l_i$-rcycles. It is
obtained from lemma \ref{C_k-conc} by doing singled out the concatenation
of fix $1$-tuples and reassembling the cycles in $S_{l_i}^{p_i}$-orbits as
in example:

$$
\begin{array}{||cc|c|c|cc|cc||}
\hline
1 & 2 & 3 & 4 & 5 & 7 & 6 & 8 \\
2 & 1 & 3 & 4 & 7 & 5 & 8 & 6 \\
\hline
\end{array}=
\begin{array}{||cc|cc|cc|cc||}
\hline
1 & 2 & 3 & 4 & 5 & 6 & 7 & 8 \\
2 & 1 & 3 & 4 & 7 & 8 & 5 & 6 \\
\hline
\end{array}\,.
$$

The difference in representations of type

$$
\begin{array}{||cc|cc||}
\hline
5 & 7 & 6 & 8 \\
7 & 5 & 8 & 6 \\
\hline
\end{array}\
\mbox{ and }\
\begin{array}{||cc|cc||}
\hline
5 & 6 & 7 & 8 \\
7 & 8 & 5 & 6 \\
\hline
\end{array}
$$
can be important if rcycles $\{\langle 57\rangle,\langle 75\rangle\}$ and
$\{\langle 68\rangle,\langle 86\rangle\}$ are not $G$-isomorphic.

There exist cases, where the partitioning of a $(l,k)$-cycle on right
$G$-related concatenation components cannot be represented as
projections of base three types $p$-orbits.  The simplest of these cases
gives the $3$-orbit of subgroup $gr((12))<S_3$.  For this case we have the
next right $G$-related $2$-orbits of $gr((12))$:

$$
\begin{array}{||cc||}
\hline
1 & 2 \\
2 & 1 \\
\hline
\end{array}\,,\
\begin{array}{||cc||}
\hline
1 & 3 \\
2 & 3 \\
\hline
\end{array}\
\mbox{ and }\
\begin{array}{||cc||}
\hline
2 & 3 \\
1 & 3 \\
\hline
\end{array}\,.
$$

The existence of such decomposition of a $(n,k)$-cycle on condition $k|n$
leeds to some intricate structures as, for example, the automorphism group
of Petersen graph (s. below).

\paragraph{Finite group permutation representation.}
Let us to consider some examples with different properties of right
permutation action. The $2$-orbit of subgroup $gr((12))<S_3$ shows an
existence of cases with no non-trivial isomorphic right permutation action
for $k$-orbits of non-normal subgroup. The example \ref{S2*S2} (s. below)
shows the existence of the right $G$-isomorphism for $2$-orbits $\{\langle
12\rangle,\langle 21\rangle\}$ and $\{\langle 34\rangle,\langle
43\rangle\}$ of normal subgroup $gr((12)(34))$ of group $S_2\otimes S_2$
that follows from proposition \ref{H->r.aut}. The next example is a
$n$-orbit of a group $G$ that is the regular permutation representation of
$S_3$ in two assemblies.

\begin{example} \label{S3(6)}
$$
\begin{array}{lcr}
\begin{array}{||cc|cc|cc||}
\hline
1 & 2 & 3 & 4 & 5 & 6 \\
\hline
2 & 1 & 5 & 6 & 3 & 4 \\
\hline
3 & 4 & 1 & 2 & 6 & 5 \\
\hline
4 & 3 & 6 & 5 & 1 & 2 \\
\hline
5 & 6 & 2 & 1 & 4 & 3 \\
\hline
6 & 5 & 4 & 3 & 2 & 1 \\
\hline
\end{array}&
\mbox{and}&
\begin{array}{||cc|cc|cc||}
\hline
1 & 2 & 3 & 5 & 4 & 6 \\
2 & 1 & 5 & 3 & 6 & 4 \\
\hline
3 & 4 & 1 & 6 & 2 & 5 \\
4 & 3 & 6 & 1 & 5 & 2 \\
\hline
5 & 6 & 2 & 4 & 1 & 3 \\
6 & 5 & 4 & 2 & 3 & 1 \\
\hline
\end{array}\,.
\end{array}
$$

\end{example}

The first table is partitioned relative to $G$-isomorphic $2$-subspaces
and the second to $S_n$-isomorphic $2$-subspaces. The example shows no
existence of a right $G$-isomorphism for $2$-rcycle $\{12,21\}$, but an
existence of right $S_6$-isomorphism for this $2$-rcycle. This fact can be
explained with next arguments: a subgroup defined by $2$-rcycle
$\{12,21\}$ has the trivial normalizer and hence its $2$-orbits of left
cosets are not necessitated to be $G$-isomorphic to the $2$-orbits of
right cosets, on the one hand, but $G$ is regular and hence necessitates
the existence of the isomorphic right action, on the other hand.

Given example shows the difference in properties of the right permutation
action in various permutation representations of a finite group. We
consider this difference and recall at first some facts from the finite
group theory.

Let $F$ be a finite group, $A<F$, $|F|/|A|=n$ and $\overrightarrow{L_n}$,
$\overrightarrow{R_n}$ be ordered partitions $FA$ and $AF$. It is known
that every transitive permutation representation of $F$ is equivalent to
the representation of $F$ given by $n$-orbits
$X_n'=\{f\overrightarrow{L_n}: f\in F\}$ or $X_n''=\{\overrightarrow{R_n}f:
f\in F\}$. It is also known that $F$ is homomorphic to its image
$Aut(X_n')$ ($Aut(X_n'')$) with the kernel of the homomorphism equal to a
maximal normal subgroup of $F$ that is contained in $A$. Further we always
assume that a finite group is isomorphic to its representation.

A maximal by inclusion subgroup $A$ of a finite group $F$ that contains no
normal subgroup of $F$ we call a \emph{md-stabilizer} of $F$ and
the corresponding representation of $F$ we call a
\emph{md-representation}. A md-stabilizer $A$ of a finite group $F$
defines a minimal degree permutation representation of $F$ in the family
of permutation representations of $F$ defined with subgroups of $A$. A
maximal degree permutation representation of $F$ is correspondingly the
permutation representation defined with trivial minimal subgroup given by
the unit of the group and this representation is called a regular
representation of $F$.

A finite group can have many (not conjugate) md-stabilizers. For example,
$S_5$ contains a transitive md-stabilizer of order $20$ generated by
permutations $(12345)$ and $(1243)$ and an intransitive md-stabilizer of
order $12$ generated by permutations $(123)(45)$ and $(23)$. The
$n$-orbits of these md-stabilizers are correspondingly:

\begin{example}\label{S_5,md}
$$
\begin{array}{lcr}
\begin{array}{||ccccc||}
\hline
1 & 2 & 3 & 4 & 5 \\
2 & 3 & 4 & 5 & 1 \\
3 & 4 & 5 & 1 & 2 \\
4 & 5 & 1 & 2 & 3 \\
5 & 1 & 2 & 3 & 4 \\
\hline
5 & 4 & 3 & 2 & 1 \\
4 & 3 & 2 & 1 & 5 \\
3 & 2 & 1 & 5 & 4 \\
2 & 1 & 5 & 4 & 3 \\
1 & 5 & 4 & 3 & 2 \\
\hline
1 & 3 & 5 & 2 & 4 \\
3 & 5 & 2 & 4 & 1 \\
5 & 2 & 4 & 1 & 3 \\
2 & 4 & 1 & 3 & 5 \\
4 & 1 & 3 & 5 & 2 \\
\hline
4 & 2 & 5 & 3 & 1 \\
2 & 5 & 3 & 1 & 4 \\
5 & 3 & 1 & 4 & 2 \\
3 & 1 & 4 & 2 & 5 \\
1 & 4 & 2 & 5 & 3 \\
\hline
\end{array}&
\mbox{ and }&
\begin{array}{||ccc|cc||}
\hline
1 & 2 & 3 & 4 & 5 \\
2 & 3 & 1 & 4 & 5 \\
3 & 1 & 2 & 4 & 5 \\
\hline
1 & 2 & 3 & 5 & 4 \\
2 & 3 & 1 & 5 & 4 \\
3 & 1 & 2 & 5 & 4 \\
\hline
1 & 3 & 2 & 4 & 5 \\
2 & 1 & 3 & 4 & 5 \\
3 & 2 & 1 & 4 & 5 \\
\hline
1 & 3 & 2 & 5 & 4 \\
2 & 1 & 3 & 5 & 4 \\
3 & 2 & 1 & 5 & 4 \\
\hline
\end{array}\,.
\end{array}
$$
\end{example}

The first md-stabilizer is a representation of group $C_5C_4=C_4C_5$. The
representation of $S_5$ with this md-stabilizer has degree $6$ equal to
maximal order of elements of $S_5$. The second md-stabilizer is a
representation of group $C_6\otimes C_2$. The representation of
$S_5$ with second md-stabilizer  (as we shall see) is the automorphism
group of Petersen graph.

Because of the property given in proposition \ref{p<n}, the special
interest is presented by md-stabilizer of a finite group with the maximal
order, which we call a \emph{least degree stabilizer} or a
\emph{ld-stabilizer} of a finite group. The corresponding representation
of a finite group we call a transitive least degree representation or
\emph{tld-representation}. That property urges also to consider a least
degree intransitive representation or a \emph{ild-representation} of a
finite group, that for some groups, for example for $C_6$, has the degree
less than the degree of tld-representation. So under a lowest degree
representation or a \emph{ld-representation} we shall understand the
smallest degree representation among tld- and ild-representations.

The given consideration puts a question: \emph{is there existing a finite
group with two non-similar ld-representations?}. If there exists no two
non-similar ld-representations, then the ld-representation is a full
invariant of a finite group and hence the study of a finite group number
invariants could be reduced to the study of ld-representation number
invariants.

For a non-minimal degree representations a simple example, of the same
degree non-similar permutation representations, is representations of $S_4$
on sets of right cosets of subgroups $gr((12))$ and $gr((12)(34))$. The
first representation contains a stabilizer of $2$-tuple on a $12$-element
set $V$ and the second contains a stabilizer of $4$-tuple on $V$.

Let $\phi:F\rightarrow S_n(V)$. We shall denote the image $\phi(F)$ of
group $F$ by representation $\phi$ as $F(V)$ and term $F(V)$ also as
representation of $F$. The $n$-orbit of $F(V)$ we term for convenience
also as representation of $F$.

\paragraph{One possible reformulation of the polycirculant conjecture.}
Let $F$ be a finite group and $A<F$ be a md-stabilizer. Let $F$ contains a
subgroup $P$ of a prime order $p$ that conjugates with no subgroup of $A$,
then it follows that $P$ is a regular subgroup of a representation $F(AF)$
and hence $p$ divides $n=|AF|$.

So the polycirculant conjecture statements that if $F(AF)$ is a $2$-closed
representation of a finite group $F$, then $F$ contains the corresponding
subgroup $P$.

To all appearance this approach cannot be successful, because it lies out
of the inside structure of a $n$-orbit.

\paragraph{Some properties of ld-representations.}
\begin{lemma}\label{A_mxB_l}
Let $A_m$ and $B_l$ be ld-representations, then the ld-representation of a
group $A_m\otimes B_l$ has degree $n=m+l$.

\end{lemma}\bigskip

The simplest example is

$$
\begin{array}{||cc|cc||}
\hline
1 & 2 & 3 & 4 \\
1 & 2 & 4 & 3 \\
\hline
2 & 1 & 3 & 4 \\
2 & 1 & 4 & 3 \\
\hline
\end{array}.
$$

\begin{theorem}\label{ild=ld}
The ld-representation of a finite group $F$ is an ild-representation if and
only if $F$ is a direct product.

\end{theorem}
\proof
If $F$ is a direct product, then the statement follows from lemma
\ref{A_mxB_l}. So let $X_n=X_{m+l}=X_m\circ X_l$ be an ild-representation
of $F$, where $X_m$ and $X_l$ are transitive, then $|X_n|/|X_m|$ and
$|X_n|/|X_l|$ greater than $1$. It follows that a stabilizer of $m$-tuple
from $X_m$ and a stabilizer of $l$-tuple from $X_l$ are normal subgroups
of $F$ (corollary \ref{fix-tuple}) and elementwise commutative.
$\Box$\bigskip

\begin{lemma}\label{ld=regular}
The ld-representation of a group $F$ is regular if and only if $F$ is not
a direct product and has a trivial $ld$-stabilizer.

\end{lemma}

\begin{corollary}\label{C_(p^m)}
The regular representation of cyclic $p$-group is ld-representation.

\end{corollary}

\begin{corollary}\label{C_n}
Let $C_n$ be ld-representation of cyclic group $C$ of order
$p_1^{m_1}\cdots p_q^{m_q}$, where $p_i$ are prime, then $C_n$ is
intransitive and $n=p_1^{m_1}+\ldots +p_q^{m_q}$.

\end{corollary}

\begin{lemma}\label{any gr.-st.}
Any finite group is a ld-stabilizer of some finite group.

\end{lemma}
\proof
Let $A$ be a finite group, then $A$ is a ld-stabilizer of a group
$F=gr(A\otimes B,d)$, where a group $B$ is isomorphic to $A$, $d$ is an
involution and $dA\otimes B=A\otimes Bd$. $\Box$\bigskip

The corresponding example gives the representation of dihedral group $D_4$:

$$
\begin{array}{||cc|cc||}
\hline
1 & 2 & 3 & 4 \\
1 & 2 & 4 & 3 \\
\hline
2 & 1 & 3 & 4 \\
2 & 1 & 4 & 3 \\
\hline
3 & 4 & 1 & 2 \\
3 & 4 & 2 & 1 \\
\hline
4 & 3 & 1 & 2 \\
4 & 3 & 2 & 1 \\
\hline
\end{array}.
$$

\paragraph{The following several sentences do the object, that we study,
more visible.}

\begin{proposition}\label{not-intersect}
Let $F$ be a finite group and $F(F)$, $F(V)$ be two its images. Let
$\alpha_k(V),\beta_k(V)$ be $k$-tuples from $V^{(k)}$ and $l$-tuples
$\alpha_l(F),\beta_l(F)$ be their images from $F^{(l)}$, then
$\alpha_k(V)$ and $\beta_k(V)$ have no intersection on $V$ if and only if
$\alpha_l(F)$ and $\beta_l(F)$ have no intersection on $F$.

\end{proposition}
\proof
Let $X_n$ be a $n$-orbit of $F(V)$ and $Y_m$ ($m=|F|$) be a $m$-orbit
of $F(F)$. The $k$-tuples $\alpha_k(V)$ and $\beta_k(V)$ are situated in
$X_n$ and $l$-tuples $\alpha_l(F)$ and $\beta_l(F)$ are situated in $Y_m$.
The statement follows from the method of reconstruction of $Y_m$ to $X_n$
that is a substitution of certain not intersected on $F$ $(m/n)$-tuples of
$Y_m$ on certain not equal elements of $V$.  $\Box$\bigskip

From this proposition follows directly

\begin{corollary}\label{Lk-Stabilizer}
Let $F$ be a finite group, $A<F$ be a md-stabilizer, $B<A$, $V=FB$ and
$G=F(V)$. Let $|G|/|A|=l$, $|A|/|B|=k$, $kl=n$, $X_n\in Orb_n(G)$,
$Y_n\subset X_n$ be a $n$-orbit of the subgroup $A(V)<G$ and
$L_n=GY_n=X_nA(V)$.  Let $I_k=AB\subset V$, $X_k=\hat{p}(I_k)X_n$ and
$Y_k=\hat{p}(I_k)Y_n$.

\begin{enumerate}
\item
Let $L_k=X_kA=GY_k=\hat{p}(I_k)L_n$, then classes of $L_k$ have no
intersection on $V$ and hence $L_k$ is a partition of $X_k$.

\item
Let $Y_n'\in L_n$, $Y_k'=\hat{p}(I_k)Y_n'$ and $Y_{n-k}'=
\hat{p}(I_{n-k})Y_n'$, then $Y_k'$ and $Y_{n-k}'$ have no intersection on
$V$.

\end{enumerate}
\end{corollary}

\begin{proposition}\label{p<n}
Let $F$ be a finite group and $p$ be a prime divisor of $|F|$. Let
$\{n_1,\,n_2,\,\ldots\}$ be degrees of transitive components of the
ld-representation of $F$, then $p\leq max(n_i)$.

\end{proposition}
\proof
The group $F$ contains a subgroup of order $p$. Any permutation $g\in S_n$
of prime order $p$ is decomposed in cycles of length either $p$ or $1$.
$\Box$\bigskip

From the definition of a minimal degree permutation representation follows

\begin{proposition}\label{reg->reg}
Let a minimal degree permutation representation contains a regular
element, then a subordinate non-minimal degree permutation representation
contains a regular element too.

\end{proposition}

\paragraph{$n$-Orbits containing $S_n$-isomorphic
$\mathbf{k}$-orbits.}

\begin{proposition}\label{Sn-iso-intr}
Let $X_n\in Orb_n(G)$, $A<G$, $Y_n\in Orb_n(A)$ and $I_k,J_k$ be
$k$-subspaces. Let $\hat{p}(I_k)Y_n$ be $S_n$-isomorphic to
$\hat{p}(J_k)Y_n$, then $\hat{p}(I_k)X_n$ is not with necessary
$S_n$-isomorphic to $\hat{p}(J_k)X_n$.

\end{proposition}
\proof
The corresponding intransitive example is simply to construct:
$$
\begin{array}{||cc|cc|cc||}
\hline
1 & 2 & 3 & 4 & 5 & 6 \\
2 & 1 & 4 & 3 & 6 & 5 \\
\hline
3 & 4 & 1 & 2 & 5 & 6 \\
4 & 3 & 2 & 1 & 6 & 5 \\
\hline
\end{array}\,.
$$
Here $2$-suborbits $\{\langle 12\rangle,\langle 21\rangle\}$ and $\{\langle
56\rangle,\langle 65\rangle\}$ are $S_n$-isomorphic, but corresponding
$2$-orbits have different power. The transitive example is not evident and
is presented in $n$-orbit of the automorphism group of Petersen graph (s.
below). $\Box$

\begin{theorem}\label{Sn-iso}
Non-minimal degree permutation representation of a finite group $F$
contains $S_n$-isomorphic $k$-orbits.

\end{theorem}
\proof
Let $B<A<F$, $|FB|=n$, $|FA|=m$ and $|AB|=k$. Let $F(FA)$ be minimal and
$F(FB)$ be non-minimal degree permutation representations of a finite
group $F$. Let $X_n$ be a $n$-orbit of $F(FB)$ and $X_m'$ be a $m$-orbit
of $F(FA)$. Let $Y_n\subset X_n$ be a $n$-orbit of $A(FB)$ and
$Y_m'\subset X_m'$ be a $m$-orbit of $A(FA)$. Let $Z_n\subset Y_n$ be a
$n$-orbit of $B(FB)$ and $Z_m'\subset Y_m'$ be a $m$-orbit of $B(FA)$. Let
$P<B$ be a subgroup of a prime order $p$. Let $T_n\subset Z_n$ be a
$n$-orbit of $P(FB)$ and $T_m'\subset Z_m'$ be a $m$-orbit of $P(FA)$.

Let $Y_k$, $Z_k\subset Y_k$ and $T_k\subset Z_k$ be $k$-orbits of $A(AB)$,
$B(AB)$ and $P(AB)$ correspondingly and let $T_{n-k}\circ T_k=T_n$.

The $m$-orbit $T_m'$, the $n$-orbit $T_n$ and the $k$-orbit $T_k$, can be
represent as a concatenation of $p$-rcycles and multituples.

A $p$-rcycle from $T_m'$ generates $k$ $p$-rcycles in $T_{n-k}$ that are
associated with cyclic permutation of $k$ right cosets of $A$. But a
multituple from $T_m'$ generates $p$-rcycles and multituple in $T_k$.
The latter $p$-rcycles are existing because $P<B<A$ and hence the action of
$P$ on $T_k$ permutes right cosets of $B<A$.

A $p$-rcycle from $T_n$, that is generated by a $p$-rcycle from $T_m'$, is
evidently not $F(FB)$-isomorphic to a $p$-rcycle from $T_n$, that is
generated by a multituple from $T_m'$. $\Box$\bigskip

This situation is demonstrated on example \ref{S3(6)}.

This property can be also emerged in a md-representation of a finite group
$F$. An example gives the group $F=C_6\otimes C_2$ in the next
representation:

\begin{example}\label{C6*C2}
$$
\begin{array}{||cc|cc|cc||}
\hline
1 & 6 & 2 & 5 & 3 & 4 \\
6 & 1 & 5 & 2 & 4 & 3 \\
\hline
2 & 1 & 3 & 6 & 4 & 5 \\
1 & 2 & 6 & 3 & 5 & 4 \\
\hline
3 & 2 & 4 & 1 & 5 & 6 \\
2 & 3 & 1 & 4 & 6 & 5 \\
\hline
4 & 5 & 3 & 6 & 2 & 1 \\
5 & 4 & 6 & 3 & 1 & 2 \\
\hline
3 & 4 & 2 & 5 & 1 & 6 \\
4 & 3 & 5 & 2 & 6 & 1 \\
\hline
5 & 6 & 4 & 1 & 3 & 2 \\
6 & 5 & 1 & 4 & 2 & 3 \\
\hline
\end{array}\,.
$$
\end{example}

It is seen that the $2$-rcycle $\{\langle 25\rangle,\langle
52\rangle\}$ does not belong to $2$-orbit of $G$, containing the
$2$-rcycle $\{\langle 16\rangle,\langle 61\rangle\}$. This matrix is a
md-representation of a finite group $F$, but not ld-representation, and it
contains a submatrix that is a non-minimal degree representation of
subgroup $S_3<F$.

In this example the $p$-subgroup does not belong to a stabilizer, but,
using the construction given in lemma \ref{any gr.-st.}, we obtain
the property for a $p$-subgroup (of order $p$) of a stabilizer of~$F$.

The considered example suggests us the next property of $n$-orbits.

\begin{proposition}\label{nmd Zn}
Let $I=\{I_k^i:i\in [1,l]\}$ be a partition of $V$ on $G$-isomorphic
$k$-subspaces and $W=\{Co(I_k): I_k\in I\}$.

Let $Z_n$ be a maximal subset of $X_n\in Orb_n(G)$ so that
$Co(\hat{p}(I_k^i)Z_n)=W$ for each $i\in [1,l]$.

Let $A<Aut(Z_n)$, $Y_n\subset Z_n$ be a $n$-orbit of $A$,
$Y_k^1=\hat{p}(I_k^1)Y_n$ be the representation $A(I_k^1)$ and
$Y_k^2=\hat{p}(I_k^2)Y_n$ be $S_{|A|}^k$-orbit, then $Z_n$ is a
non-minimal degree representation.

\end{proposition}
\proof
$Z_n$ is automorphic, because $GZ_n$ is a partition of $X_n$. From
proposition \ref{H->r.aut} follows that $A$ contains no normal subgroup of
$Aut(Z_n)$. Hence $Z_n$, that is defined on $V$, is isomorphic to a
$l$-orbit $Z_l'$, that is defined on $W$ and obtained by the evident
reduction of $Z_n$. $\Box$\bigskip

The next example of a minimal degree representation of the group
$S_5\otimes S_2$ contains $S_n$-isomorphic $p$-orbits in the case
$(p,n)=1$.

\begin{example}\label{S5*S2}
$$
\begin{array}{||c|cc|cc||}
\hline
1 & 2 & 5 & 3 & 4 \\
1 & 5 & 2 & 4 & 3 \\
\hline
2 & 3 & 1 & 4 & 5 \\
2 & 1 & 3 & 5 & 4 \\
\hline
3 & 4 & 2 & 5 & 1 \\
3 & 2 & 4 & 1 & 5 \\
\hline
4 & 5 & 3 & 1 & 2 \\
4 & 3 & 5 & 2 & 1 \\
\hline
5 & 1 & 4 & 2 & 3 \\
5 & 4 & 1 & 3 & 2 \\
\hline
\end{array}\,.
$$
\end{example}

This example shows an existence of $n$-orbit of a subgroup of order $p$,
that contains $S_n$-isomorphic $p$-rcycles, but no $S_p^p$-orbit
$G$-isomorphic to a $p$-rcycle. Here: $2$-rcycles $\{\langle
25\rangle,\langle 52\rangle\}$ and $\{\langle 34\rangle,\langle
43\rangle\}$ are $S_n$-isomorphic, $S_p^p$-orbits $\langle\langle
23\rangle\langle 54\rangle\rangle$ and $\langle\langle 54\rangle\langle
23\rangle\rangle$ are equal and hence $G$-isomorphic and $2$-rcycle
$\{\langle 25\rangle,\langle 52\rangle\}$ is $G$-isomorphic to a
$2$-orbit $\{\langle 13\rangle,\langle 14\rangle\}$. We shall see that
properties of this $5$-orbit give an appearance to unconventional
properties of the $10$-orbit of the Petersen graph automorphism group.

\begin{theorem}\label{(p,n)=1,mdr}
Let a prime $p$ divides $|G|$ and does not divide $n$, then $G$ is a
md-representation.

\end{theorem}
\proof
Let $(\alpha_p)$ be a cycle, $X_p=G\alpha_p$, $L_p=G(\alpha_p)\alpha_p$
and $L_1=\hat{p}(I_1)L_p$, then $L_1=Co(L_p)$ is a covering of $V$ on
$G$-isomorphic $p$-subsets. We consider two possible cases.

\begin{enumerate}
\item
Let $\sqcup L_1=V$, then $G$ is a primitive group and hence there exists no
partition $Q$ of $V$ on $G$-isomorphic subsets for that $G(Q)$ would be a
representation of $G(V)$. Indeed, if $g\in G$ is a permutation of order
$p$, then $gQ\sqcap Q\neq Q$ for any partition $Q$. It contradicts
proposition~\ref{not-intersect}.

\item
Let now $\sqcup L_1=Q$, where $Q$ is a partition of $V$, then $Q$ consists
of $G$-isomorphic classes. But in this case, because of transitivity $G$
and $(n,p)=1$, there exists a $n$-orbit of a $p$-subgroup whose projections
on subspaces from $Q$ are $G$-isomorphic and hence $X_n$ does not contain
$S_n$-isomorphic $(n/|Q|)$-orbits. Thus, accordingly to theorem
\ref{Sn-iso}, $X_n$ is a minimal degree permutation representation. $\Box$

\end{enumerate}\bigskip

An example of the second case representation of a finite group $F$ can be
obtained from lemma~\ref{any gr.-st.}, if to assign $A=A_4$ and $p=3$. The
consideration of this example for $A=S_3$ and $p=2$ shows that the
condition, $p$ does not divide $n$, is not of principal for the second case
of theorem~\ref{(p,n)=1,mdr}. We shall see that namely this situation
takes place in the $10$-orbit of the Petersen graph automorphism group.

\paragraph{Conditions of $k$-closure and properties of $k$-closed groups.}

\begin{proposition}\label{Xk,Yk Aut}
Let $Y_k\in Orb_k(Aut(X_k))$, then it is not follows that
$Aut(Y_k)=Aut(X_k)$.

\end{proposition}
\proof
An example: $X_2=\{\langle 12\rangle,\langle 23\rangle,\langle
34\rangle,\langle 41\rangle\}$, $Y_2=\{\langle 13\rangle,\langle
31\rangle,\langle 24\rangle,\langle 42\rangle\}$.
$\Box$

\begin{proposition}\label{Xk,Yk iso}
Let $k$-orbits $Y_k$ and $X_k$ be isomorphic and $Aut(X_k)=Aut(Y_k)$, then
it is not follows that $Y_k=X_k$.

\end{proposition}
\proof
An example: $X_2=\{\langle 13\rangle,\langle 24\rangle\}$, $Y_2=\{\langle
14\rangle,\langle 23\rangle\}$. $\Box$

\begin{proposition}\label{iso-rorbits}
Let $X_k$ and $Y_k$ be isomorphic $k$-rorbits with the same automorphism
group, then it is not follows that $X_k=Y_k$.

\end{proposition}
\proof
The $2$-orbits $X_2=\hat{p}(\langle 23\rangle)X_5$ and $Y_2=\hat{p}(\langle
45\rangle)X_5$ from example \ref{S5*S2} represent such case.
$\Box$

\begin{lemma}\label{k.cl.gr->k.cl.subgr}
Let $X_n$ be a $n$-orbit of a $k$-closed group $G$, $A<G$ and $Y_n\subset
X_n$ be a $n$-orbit of $A$. Let $Y_k(I_k)=\hat{p}(I_k)Y_n$,
$B=\cap_{(I_k\subset I_n)}Aut(Y_k(I_k))$, $P_n=G\,(\cup BY_n)$ and classes
of $P_n$ have no intersections, then $A$ is $k$-closed.

\end{lemma}
\proof
Let $X_k(I_k)=\hat{p}(I_k)X_n$. It is given that $G=\cap_{(I_k\subset
I_n)}Aut(X_k(I_k))$. Further we have: $Aut(Y_n)<B$, $Y_n\in BY_n$,
$Y_n\subset\cup BY_n$ and $X_n=\cup GY_n\subseteq \cup G\cup BY_n=\cup
P_n$.  Let $Y_n$ be not $k$-closed, then every class of $GY_n$ is not
$k$-closed.  As classes of $P_n$ have no intersections, then $X_n\subset
\cup P_n$ and hence $X_n$ is not $k$-closed.  Contradiction.  $\Box$

\begin{theorem}\label{Sn.iso->2-cl}
Let a transitive $n$-orbit $X_n$ contains $S_n$-isomorphic $k$-projections
($k\geq 2$), then $X_n$ is $2$-closed.

\end{theorem}
\proof
We have two possibilities:

\begin{enumerate}
\item
There exists a subspace $I_4=\langle 1234\rangle$ so that $4$-orbit
$X(1234)=\hat{p}(\langle 1234\rangle)X_n$ is not $2$-closed and $2$-orbits
$X(12)=\hat{p}(\langle 12\rangle)X_n$ and $X(34)=\hat{p}(\langle
34\rangle)X_n$ are $S_n$-isomorphic. Then $3$-orbits
$X(123)=\hat{p}(\langle 123\rangle)X(1234)$ and $X(234)=\hat{p}(\langle
234\rangle)X(1234)$ are also not $2$-closed.

Let $A(ijk\ldots)=Aut(X(ij))\cap Aut(X(ik))\cap Aut(X(jk))\cap\,\ldots\,$,
then $\cup Aut(123)X(1234)$ and $\cup Aut(234)X(1234)$ have to be
$2$-closed and equal, i.e. $Aut(123)=Aut(234)=Aut(1234)$. But such
equality (for transitive $4$-orbit $X(1234)$) is impossible, because
$Aut(123)$ and $Aut(234)$ are conjugate subgroups of $S_n$ and hence are
not equal.

\item
There exists a subspace $I_6=\langle 123456\rangle$ so that $6$-orbit
$X(123456)=\hat{p}(\langle 123456\rangle)X_n$ is not $2$-closed,
$3$-orbits $X(123)=\hat{p}(\langle 123\rangle)X_n$ and
$X(456)=\hat{p}(\langle 456\rangle)X_n$ are $S_n$-isomorphic and not
$2$-closed. Then $\cup Aut(123)X(123456)$ and $\cup Aut(456)X(123456)$
have to be $2$-closed and equal or $Aut(123)=Aut(456)=Aut(123456)$. This
equality is also impossible for the same reason. $\Box$

\end{enumerate}

Let $X_n$, $Y_m$ and $Z_l$ be three representations of a finite group $F$.
Let $n>m>l$. It is of interest the relation between a $k$-closure
property of $Z_l$, $Y_m$ and $X_n$.

It is evident that, if $Y_m$ is not $1$-closed, then $X_n$ is not
$1$-closed too and $Z_l$ can be $1$-closed. And, if $Y_m$ is $k$-closed for
$k>1$, then $X_n$ is $k$-closed too and $Z_l$ can be not $k$-closed.

\paragraph{Unconventional cyclic structure on $k$-orbits.}
Now we shall consider one interesting property of the right permutation
action on $k$-orbits that has no analogy in the group theory. The right
automorphism action on a $n$-orbit $X_n$ of a group $G$ maps a
$k$-subspace $I_k$ on an isomorphic $k$-subspace $I_k'$ and so
$X_k=\hat{p}(I_k)X_n=\hat{p}(I_k')X_n=X_k'$.  The latter equality
generates an unconventional cyclic structure on a $k$-orbit $X_k$, that
one can see on next examples:

\begin{example}\label{S2*S2}
$$
\left(
\begin{array}{cc|cc}
1 & 2 & 3 & 4 \\
2 & 1 & 4 & 3 \\
\hline
3 & 4 & 1 & 2 \\
4 & 3 & 2 & 1 \\
\end{array}
\right)
\rightarrow
\left(
\begin{array}{c|c}
B & B' \\
\hline
B'& B \\
\end{array}
\right)
\mbox{ \emph{and}  }
\left(
\begin{array}{ccc}
1 & 2 & 3 \\
2 & 1 & 3 \\
1 & 3 & 2 \\
3 & 1 & 2 \\
2 & 3 & 1 \\
3 & 2 & 1 \\
\end{array}
\right)
\rightarrow
\left(
\begin{array}{cc|cc}
1 & 2 & 2 & 3 \\
2 & 1 & 1 & 3 \\
\hline
2 & 3 & 3 & 1 \\
1 & 3 & 3 & 2 \\
\hline
3 & 1 & 1 & 2 \\
3 & 2 & 2 & 1 \\
\end{array}
\right)
\rightarrow
\left(
\begin{array}{c|c}
B_1 & B_2 \\
\hline
B_2 & B_3 \\
\hline
B_3 & B_1 \\
\end{array}
\right).
$$

\end{example}

We have introduced the right permutation action on $n$-orbits as a
permutation of coordinates of $n$-tuples. Of course, we can consider this
action as the permutation $n$-tuples with the same result. Such
interpretation of the right permutation action leads to next

\begin{lemma}\label{Rk=Yn}
Let $X_k$ be a $k$-orbit of a group $G$, $A<G$, $R_k=AX_k$ be a
partition of $X_k$ on $k$-orbits of right cosets of $A$ in $G$ and $Y_n$
be $n$-orbit of $A$. Then for any $k$-orbit $Y_k\in R_k$ there exists a
subspace $I_k$ so that $Y_k=\hat{p}(I_k)Y_n$.

\end{lemma}
\proof
It follows from $A\alpha_n=Y_n$ for $\alpha_n\in Y_n$ and
$A\alpha_k\in R_k$ for $\alpha_k\in X_k$.  $\Box$\bigskip

\begin{theorem}\label{rr}
Let $X_n$ be a $n$-orbit of a group $G$, $I_k,I_k'$ be isomorphic
subspaces, and $X_{2k}=\hat{p}(I_k\circ I_k')X_n=\uplus X_k\circ \uplus
X_k$. Let $A<G$, $R_{2k}=AX_{2k}$ and $R_k=AX_k$, then $R_{2k}$ can be
partition in cycles on classes of $R_k$.

\end{theorem}
\proof
It follows from lemmas \ref{Rk=Yn} and \ref{mMmM}. $\Box$\bigskip

Namely this property we can see on examples.

\begin{remark}\label{rem}
Let $I_k^1$, $I_k^2$ and $I_k^3$ be isomorphic subspaces of a $n$-orbit
$X_n$, then $\hat{p}(I_k^1)X_n=\hat{p}(I_k^2)X_n=\hat{p}(I_k^3)X_n$. But
on the condition it does not follow that $\hat{p}(I_k^1\circ
I_k^2)X_n=\hat{p}(I_k^2\circ I_k^3)X_n$.

Hence the corresponding cycle structure on $lk$-orbits of right cosets of
a subgroup does not exist with necessary for $l>2$. This fact represents
the difference between $2$-closed groups and $m$-closed groups for $m>2$.

\end{remark}

\begin{proposition}\label{rlrl}
Let $X_{2k}\in Orb_{2k}(G)$, $A<G$, $R_{2k}=AX_{2k}$,  $Y_k^1$ and
$Y_k^2$ be isomorphic $k$-orbits of a subgroup $A$, $Y_{2k}=Y_k^1\circ
Y_k^2\in R_{2k}$ and $C_{2k}\subset R_{2k}$ be a cycle containing
$Y_{2k}$. Let the set $Z_{2k}=\cup C_{2k}$ is a $(2k)$-orbit of some
subgroup $B<G$, then $A$ is a normal subgroup of $B$.

\end{proposition}
\proof
Since $Y_k^1$ and $Y_k^2$ are isomorphic, $k$-orbits $Y_k^i$, that form
the cycle $C_{2k}$, are $k$-orbits of left and right cosets of $A$ in $B$.
Then statement follows from lemma \ref{Lk*Rk}. $\Box$


\section{Correspondence between $k$-orbits and their automorphism\\ groups}

\begin{lemma}
Every cycle $c\subset\in G$ of length $l\geq k$ corresponds to a
$(l,k)$-cycle $C_{lk}$ of some $k$-orbit $X_k\in Orb_k(G)$.

\end{lemma}
\proof
If $c=(\alpha_l)$, then $l$-orbit $X_l=G\alpha_l$ contains $l$-rcycle
$rC_l=gr(c)\alpha_l$.  The $X_k$ and $C_{lk}$ are corresponding
$k$-projections of $X_l$ and $rC_l$. $\Box$\bigskip

For $k>l$ the counterexample is given by the cycle $(56)$ in fourth case
of example \ref{simple-struc}, where there exists no $(2,3)$-cycle for
subautomorphism $(56)$.

The reverse statement for $k$-orbits of not $k$-closed groups is not
correct. An example is $2$-orbit of $A_4$ that contains $(4,2)$-cycle
related to no subautomorphism of $A_4$.

For $k$-closed groups the reverse statement is also not correct. This
shows an example of $2$-closed group that is defined by $2$-orbit
$X_2=\{14,25,36,41,52,63\}$. The group $Aut(X_2)$ has a $2$-orbit
$X_2'=\{12,13,15,16,21,24,23,26,32,35,31,34,42,45,43,46,51,54,53,
56,62,65,61,64\}$.\label{X_2'} The automorphism group of $2$-suborbit
$\{12,13,15,16\}\subset X_2'$ contains a cycle $(2536)$ that does not
belong to $Aut(X_2')$. We can see that the possibility for construction of
this counterexample gives a concatenation of $S_4^1$-orbit with
$(4,1)$-multituple. But $X_2'$ contains a suborbit $\{12,24,43,31\}$ that
is a projection of a $4$-rcycle and also is not a $2$-orbit of a subgroup
of $Aut(X_2')$. In latter case the length $l$ of a cycle $(1243)$ is not
prime first and does not divide degree $n$ second. For a prime $l$ that
does not divide $|Aut(X_2')|=48$ we have the next counterexample: an
automorphic $2$-subset $\{13,32,24,45,51\}\subset X_2'$.


\subsection{The local property of $k$-orbits}
The trivial case of reverse statement we obtain from
corollary \ref{Y_k-A}. For disclosing of a non-trivial local property of
an automorphism group $k$-orbits we have to consider the case, where
$(p,k)$-subset of a $k$-orbit $X_k$ is a $k$-projection of $p$-rcycle for
$p$ being a prime divisor of $|Aut(X_k)|$. We shall write further a
$k$-projection of $l$-rcycle for $k\leq l$ as \emph{$(l,k)$-rcycle}
$rC_{lk}$.

\begin{theorem}\label{rC_pk->A}
Let $X_k$ be automorphic, $p\geq k$ be a prime divisor of $|Aut(X_k)|$ and
$rC_{pk}\subset X_k$ be a $(p,k)$-rcycle, then $Aut(rC_{pk})<<Aut(X_k)$.

\end{theorem}
\proof
Let $L_k=Aut(X_k)rC_{pk}$ and $|L_k|p$ divides $|Aut(X_k)|$, then the
statement follows from theorem \ref{|Y_k||GY_k|}. Let $|L_k|p>|Aut(X_k)|$,
then $|L_k|=|Aut(X_k)|$ and hence $L_k$ can be partition on subsets
$L_k^i$, $i\in [1,p]$ so that $|L_k^i|p=|Aut(X_k)|$. Since subsets $L_k^i$
have no intersections, there exist subgroups $A_i<Aut(X_k)$ so that
$A_iL_k^i=L_k^i$. It follows that $|L_k|p$ divides $|Aut(X_k)|$ and hence
$|L_k|p=|Aut(X_k)|$.  $\Box$\bigskip

The theorem \ref{rC_pk->A} and lemma \ref{XAUXB=XAB} give a
possibility for the reconstruction of subautomorphisms of $k$-orbit
through its symmetry properties.

The next statement gives the relation between automorphism group of a
$k$-orbit and automorphism group of its $k$-suborbit.

\begin{proposition}\label{Aut(Y_k)<<Aut(X_k)}
Let $X_k$ be a $k$-orbit and $Y_k\subset X_k$ be a $k$-orbit of a subgroup
$A=Aut(Y_k;V)\cap Aut(X_k)$, then $Aut(Y_k)<<A<Aut(X_k)$ if and only if
$Aut(Y_k)<<\cap_{(Z_k\in AX_k)}Aut(Z_k;V)$.

\end{proposition}
\proof
The statement follows from evident equality $A=\cap_{(Z_k\in
AX_k)}Aut(Z_k;V)$. $\Box$\bigskip


\section{Primitivity and imprimitivity.}

Below a group $G$ is transitive.

Let $X_k$ be a $k$-rorbit, $Y_k\subset X_k$ be a $k$-subrorbit,
$\alpha_k\in Y_k$ and $Co(Y_k)=Co(\alpha_k)$, then $L_k=GY_k$ is a
partition of $X_k$ (lemma~\ref{Co(Y_k)=Co(alpha_k)}), but classes of $L_k$
can be intersected on $V$. We shall call a $k$-rorbit $X_k$ for $k<n$
\emph{$V$-coherent}, if $\sqcup Co(X_k)=V$ and \emph{$V$-incoherent} if
$\sqcup Co(X_k)$ is a partition of $V$. We shall write simply coherent and
incoherent, instead of $V$-coherent and $V$-incoherent, if it will be
clear, what a set we consider.

\begin{proposition}\label{Aut-incoherent}
The automorphism group of an incoherent $k$-rorbit is imprimitive.

\end{proposition}

\begin{corollary}\label{impr-gr}
A group $G$ is imprimitive if and only if it contains an incoherent
$k$-rorbit.

\end{corollary}

\begin{corollary}\label{nmd-imprimit}
Non-minimal degree representations are imprimitive.

\end{corollary}

The automorphism group of a coherent $k$-rorbit can be imprimitive. The
example is the $2$-rorbit $\{13,31,24,42,14,41,23,32\}$.

Let a coherent (incoherent) $k$-rorbit $X_k$ contains no $V$-coherent and
no $V$-incoherent $k$-subrorbit, then $X_k$ can be called an
\emph{elementary coherent} (\emph{elementary incoherent}) $k$-rorbit.

\begin{lemma}\label{el.coherent-primitive}
The automorphism group of an elementary coherent $k$-rorbit is primitive.

\end{lemma}

\begin{corollary}\label{primitive-el.coherent}
The group $G$ is primitive if and only if it contains an elementary
$V$-coherent $k$-subrorbit.

\end{corollary}

A maximal $k$-subrorbit $Y_k$ of a $k$-rorbit $X_k$, that is a structure
element of coherent (incoherent) $k$-subrorbits, we call a
\emph{$k$-block}.

Let $Y_k$ be a $k$-block of an incoherent $k$-rorbit $X_k$, then
$U=Co(Y_k)$ is a $1$-block or $k$-element block of an imprimitive group
$G$ in conventional definition.

Let us to give some examples of coherent and incoherent $k$-rorbits:

\begin{enumerate}
\item
A $2$-orbit $X_2$ of $S_3$ is elementary coherent and a $2$-rcycle from
$X_2$ is a $2$-block.

\item
$2$-orbit of $C_5\otimes C_2$ is elementary coherent. This group contains
$S_5$-isomorphic elementary coherent $2$-rorbits.

\item
A $2$-orbit of $A_5$ is coherent but not elementary coherent. It
contains an elementary coherent $2$-orbit of $C_5\otimes C_2$. A
$3$-orbit of $A_5$ is elementary coherent and contains elementary coherent
suborbits on $4$-element subsets of $V$.

\item
All $2$-orbits of $C_2\otimes C_2$ and two from six $2$-orbits of $D_4$
are elementary incoherent. Other four $2$-orbits of $D_4$ are coherent.

\end{enumerate}

\begin{proposition}\label{el.impr->base.typr}
Let $p$ be prime and $X_p$ be an elementary incoherent $p$-rorbit, then
there exists a subgroup $A<Aut(X_p)$ of order $p$ for that classes of a
partition $R_p=AX_p$ are $p$-orbits of the base type, i.e. they are either
$p$-rcycles, or $S_p^p$-orbits or $p$-tuples.

\end{proposition}
\proof
An elementary incoherent $p$-rorbit consists of not intersected on $V$
$p$-rcycles. $\Box$\bigskip

The reverse statement:

\begin{lemma}\label{base.typ->impr}
Let $p$ be prime, $X_p$ be a $p$-rorbit and $rC_p\subset X_p$ be a
$p$-rcycle that is a $p$-orbit of a subgroup $A<Aut(X_p)$ of order $p$.
Let classes of $R_p=AX_p$ be $p$-orbits of the base type, then $Aut(X_p)$
is imprimitive.

\end{lemma}
\proof
Let the statement is not correct and $Aut(X_p)$ be primitive, then
$L_p=Aut(X_p)rC_p$ contains intersected on $V$ classes and hence there
exists a class of $R_p$ that has a $(p,m<p)$-multituple as a concatenation
component. Contradiction. $\Box$\bigskip

In the next example:
$L_1=\{\{124\},\{235\},\{346\},\{451\},\{562\},\{613\}\}$, we see that
$\sqcup L_1=V=[1,6]$, where $k=3$ divides $n=6$, but the corresponding
$3$-set $X_3=\{124,241,421,\ldots\}$ is not automorphic.

\begin{theorem}\label{not k|n}
Let $X_k$ be an elementary coherent $k$-rorbit, then $\neg (k|n)$.

\end{theorem}
\proof
The statement is correct for $n$ being a prime, because of $k<n$, so we
assume that $n$ is not prime. Let the statement is not correct for some
$k$, then it is not correct also for a prime divisor $p$ of $k$. So we
assume that $k=p$ is a prime. Let $Y_p\subset X_p$ be a $p$-rcycle that is
a $p$-orbit of a subgroup $A<Aut(X_p)$ of order $p$, $R_p=AX_p$ and
$L_p=Aut(X_p)Y_p$. Since by hypothesis classes of $L_p$ are intersected on
$V$, then $R_p$ contains classes that have a $(p,l<p)$-multituple as a
concatenation component. Let $X_n\in Orb_n(Aut(X_p))$ and $Y_n\subset X_n$
be a $n$-orbit of $A$, then $Y_n$ is a concatenation of (not with
necessary $Aut(X_p)$-isomorphic) $p$-rcycles and a $(p,pr)$-multituple.
Let $I_{pr}$ be a subspace defining $(p,pr)$-multituple, then, accordingly
to lemma \ref{k-rorbit}, the $pr$-orbit $X_{pr}=\hat{p}(I_{pr})X_n$ is a
$pr$-rorbit and hence $l|pr$. Another, the subspace $I_{pr}$ is a
concatenation of $m=pr/l$ $Aut(X_p)$-isomorphic subspaces $I_l^i$ ($i\in
[1,m]$). It follows that $R_p$ contains $m$ $Aut(X_p)$-isomorphic classes
and that $Aut(X_p)$ contains a subgroup $B$ that acts transitive on the
system of these $m$ classes. It follows that $R_p$ and hence $X_p$ can be
partition on $m$ $Aut(X_p)$-isomorphic classes. Let $L_p'=L_p\sqcup R_p$
be this partition of $X_p$ on $m$ $Aut(X_p)$-isomorphic classes, then
classes of $L_p'$ are $p$-orbits of a normal subgroup of $Aut(X_p)$ and
cannot be intersected on $V$. Hence $X_p$ is not elementary coherent
$p$-rorbit.  Contradiction. $\Box$


\section{Petersen graph}

Here we consider the properties of $n$-orbits that are not visible
from the group theory and therefore had hindered to solve the
polycirculant conjecture. It is the cases, where $p|n$, but there exists no
transitive and imprimitive subgroup of order $p$ and therefore there
exists no subgroup of order $p$, whose $n$-orbits could be represented as
a concatenation of $(p,p)$-orbits of base type.

The simplest example is a $(n=6)$-orbit of a group $G=S_3\otimes C_2$.

\begin{example}\label{S3xC2}
$$
\begin{array}{||cc|c|cc|c||}
\hline
1 & 2 & 3 & 4 & 5 & 6 \\
2 & 1 & 3 & 5 & 4 & 6 \\
\hline
1 & 3 & 2 & 4 & 5 & 6 \\
3 & 1 & 2 & 5 & 4 & 6 \\
\hline
2 & 3 & 1 & 5 & 6 & 4 \\
3 & 2 & 1 & 6 & 5 & 4 \\
\hline
4 & 5 & 6 & 1 & 2 & 3 \\
5 & 4 & 6 & 2 & 1 & 3 \\
\hline
4 & 5 & 6 & 1 & 3 & 2 \\
5 & 4 & 6 & 3 & 1 & 2 \\
\hline
5 & 6 & 4 & 2 & 3 & 1 \\
6 & 5 & 4 & 3 & 2 & 1 \\
\hline
\end{array}\,.
$$

\end{example}

It is seen that the pair $\langle 12\rangle$ is $G$-isomorphic to $\langle
45\rangle$, but not $G$-isomorphic to $\langle 36\rangle$, so $6$-orbit of
a subgroup $gr((12)(45))$ cannot be represented as a concatenation of
$p$-orbits of base type. We shall say that this subgroup of a prime order
and its $n$-orbit are \emph{undecomposable}. Nevertheless the given group
contains a \emph{decomposable} (on the $(p,p)$-orbits of the base type)
subgroup $gr((14)(25)(36))$, where $\{14\},\{25\},\{36\}$ are incoherent
$2$-blocks of a corresponding imprimitive subgroup of $G$.

The Petersen graph gives an example, where for the least prime divisor
$p=2$ of the degree $n=10$ there exists no decomposable subgroup with
incoherent $2$-blocks. From here follows unconventional properties of the
automorphism group of this graph. The automorphism group of Petersen graph
is a representation $G$ of $S_5$ on $10$-element set $V$. It is a
representation of $S_5$ with right (left) cosets of a subgroup of order
$12$ represented in example \ref{S_5,md}. This representation can be also
obtained by action of $S_5$ on unordered pairs from the set
$V'=\{1,2,3,4,5\}$, where $V=\{1=\{1,2\},2=\{1,3\},\ldots,0=\{4,5\}\,\}$.

The following matrix is a $10$-orbit $Y_{10}$ of a transitive, imprimitive
subgroup of $G$ (that is isomorphic to a first subgroup from example
\ref{S_5,md}).

\begin{example}
$$
\begin{array}{||ccccc|ccccc||c||}
\hline
12 & 23 & 34 & 45 & 15   & 14 & 24 & 25 & 35 & 13 & \\
\hline
\hline
1  & 5  & 8  & 0  & 4    & 3  & 6  & 7  & 9  & 2     & e \\
5  & 8  & 0  & 4  & 1    & 7  & 9  & 2  & 3  & 6     & (12345) \\
8  & 0  & 4  & 1  & 5    & 2  & 3  & 6  & 7  & 9     & (13524) \\
0  & 4  & 1  & 5  & 8    & 6  & 7  & 9  & 2  & 3     & (14253) \\
4  & 1  & 5  & 8  & 0    & 9  & 2  & 3  & 6  & 7     & (15432) \\
\hline
3  & 6  & 7  & 9  & 2    & 4  & 0  & 8  & 5  & 1      & (2453) \\
6  & 7  & 9  & 2  & 3    & 8  & 5  & 1  & 4  & 0      & (1435) \\
7  & 9  & 2  & 3  & 6    & 1  & 4  & 0  & 8  & 5      & (1254) \\
9  & 2  & 3  & 6  & 7    & 0  & 8  & 5  & 1  & 4      & (1523) \\
2  & 3  & 6  & 7  & 9    & 5  & 1  & 4  & 0  & 8      & (1342) \\
\hline
4  & 0  & 8  & 5  & 1    & 2  & 9  & 7  & 6  & 3      & (25)(34) \\
0  & 8  & 5  & 1  & 4    & 7  & 6  & 3  & 2  & 9      & (15)(24) \\
8  & 5  & 1  & 4  & 0    & 3  & 2  & 9  & 7  & 6      & (23)(14) \\
5  & 1  & 4  & 0  & 8    & 9  & 7  & 6  & 3  & 2      & (45)(13) \\
1  & 4  & 0  & 8  & 5    & 6  & 3  & 2  & 9  & 7      & (12)(35) \\
\hline
2  & 9  & 7  & 6  & 3    & 1  & 5  & 8  & 0  & 4      & (2354) \\
9  & 7  & 6  & 3  & 2    & 8  & 0  & 4  & 1  & 5      & (1325) \\
7  & 6  & 3  & 2  & 9    & 4  & 1  & 5  & 8  & 0      & (1534) \\
6  & 3  & 2  & 9  & 7    & 5  & 8  & 0  & 4  & 1      & (1243) \\
3  & 2  & 9  & 7  & 6    & 0  & 4  & 1  & 5  & 8      & (1452) \\
\hline
\end{array}
$$
\end{example}
and the reassembling of this example:
\begin{example}
$$
\begin{array}{||cc|cc|cc|cc|cc||c||}
\hline
12 & 15 & 14 & 13 & 23 & 35 & 34 & 45 & 24 & 25  & \\
\hline
\hline
1  & 4  & 3  & 2  & 5  & 9  & 8  & 0  & 6  & 7      & e \\
4  & 1  & 2  & 3  & 0  & 6  & 8  & 5  & 9  & 7      & (25)(34) \\
\hline
5  & 1  & 7  & 6  & 8  & 3  & 0  & 4  & 9  & 2      & (12345) \\
1  & 5  & 6  & 7  & 4  & 9  & 0  & 8  & 3  & 2      & (12)(35) \\
\hline
8  & 5  & 2  & 9  & 0  & 7  & 4  & 1  & 3  & 6      & (13524) \\
5  & 8  & 9  & 2  & 1  & 3  & 4  & 0  & 7  & 6      & (45)(13) \\
\hline
0  & 8  & 6  & 3  & 4  & 2  & 1  & 5  & 7  & 9      & (14253) \\
8  & 0  & 3  & 6  & 5  & 7  & 1  & 4  & 2  & 9      & (23)(14) \\
\hline
4  & 0  & 9  & 7  & 1  & 6  & 5  & 8  & 2  & 3      & (15432) \\
0  & 4  & 7  & 9  & 8  & 2  & 5  & 1  & 6  & 3      & (15)(24) \\
\hline
3  & 2  & 4  & 1  & 6  & 5  & 7  & 9  & 0  & 8      & (2453) \\
2  & 3  & 1  & 4  & 9  & 0  & 7  & 6  & 5  & 8      & (2354) \\
\hline
6  & 3  & 8  & 0  & 7  & 4  & 9  & 2  & 5  & 1      & (1435) \\
3  & 6  & 0  & 8  & 2  & 5  & 9  & 7  & 4  & 1      & (1452) \\
\hline
7  & 6  & 1  & 5  & 9  & 8  & 2  & 3  & 4  & 0      & (1254) \\
6  & 7  & 5  & 1  & 3  & 4  & 2  & 9  & 8  & 0      & (1243) \\
\hline
9  & 7  & 0  & 4  & 2  & 1  & 3  & 6  & 8  & 5      & (1523) \\
7  & 9  & 4  & 0  & 6  & 8  & 3  & 2  & 1  & 5      & (1534) \\
\hline
2  & 9  & 5  & 8  & 3  & 0  & 6  & 7  & 1  & 4       & (1342) \\
9  & 2  & 8  & 5  & 7  & 1  & 6  & 3  & 0  & 4       & (1325) \\
\hline
\end{array}
$$
\end{example}

It can be seen that there exists no partition of $2$-projection
$\hat{p}(\langle 14\rangle)Y_{10}$ on not intersected on $V$ classes, but
this covering of $V$ with $2$-tuples can be partition in two not
intersected coverings:  $\{14,15,58,40,08\}$ and $\{23,36,29,79,67\}$ that
form elementary coherent $2$-orbits on corresponding two $5$-element
subsets of $V$. This case is similar to we could see on example
\ref{S5*S2}.

The following example of a $(2,10)$-orbit

$$
\begin{array}{||cc||cc|cc|cc|cc||}
\hline 24 & 35 &
12 & 14 & 13 & 15 & 23 & 25 & 45 & 34 \\
\hline
\hline
6  & 9  & 1  & 3  & 2  & 4  & 5  & 7  & 0  & 8  \\
9  & 6  & 3  & 1  & 4  & 2  & 7  & 5  & 8  & 0  \\
\hline
\end{array}
$$
contains $S_{10}$ isomorphic $2$-rcycles
$\{\langle 69\rangle,\langle 96\rangle\}$ and $\{\langle 13\rangle,\langle
31\rangle\}$, but the corresponding two projections of $X_{10}$ are not
$S_{10}$ isomorphic, because the $2$-orbit $X_2(\langle 69\rangle)$
consists of $30$ pairs and $2$-orbit $X_2(\langle 13\rangle)$ consists of
$60$ pairs. It is a transitive realization of the property of example from
proposition \ref{Sn-iso-intr}.

\begin{remark}
One can see that given properties of $X_n$ cannot be obtained in the group
theory, because they are properties of the internal structure of $X_n$. Of
course, the internal structure of $X_n$ characterizes the group and
therefore its properties are also group properties. But these properties
of a group lie out the group algebra that characterizes $X_n$ as whole.

\end{remark}


\section{The proof of the polycirculant conjecture}

\subsection{The proof of lemma \ref{tr.impr<tr.pr}}

The proof of lemma \ref{tr.impr<tr.pr} follows from

\begin{lemma}\label{qGDn}
Let $q$ be the greatest prime divisor of $n$, then $G$ contains a
transitive, imprimitive subgroup with incoherent $q$-blocks.

\end{lemma}
\proof
Let the statement is not correct, then every $q$-rorbit $X_q(I_q)$
contains a coherent $q$-subrorbit $Y_q(I_q)$ on some automorphic
$k$-subspace $I_k$, where $k$ is a divisor  of $n$ greater than $q$,
$I_q\subset I_k$ and $X_k(I_k)$ is an incoherent $k$-rorbit. Let $p$ be a
prime divisor of $k$, then it follows that $p<q$ and hence there exists an
automorphic $p$-subspace $I_p\subset I_q$. Therefore there exists a
coherent $p$-subrorbit $Y_p(I_p)$ of a $p$-rorbit $X_p(I_p)$ on the
$q$-subspace $I_q$. Since $q|n$, classes of the partition $L_p=GY_p(I_p)$
of $X_p(I_p)$ are not intersected on $V$ and hence $X_q(I_q)$ can not
contain a coherent $q$-subrorbit $Y_q(I_q)$ on the $k$-subspace $I_k$.
Contradiction.  $\Box$

\subsection{The proof of theorem \ref{2cl.impr->reg}}

Let $G$ be a transitive group, then it contains a transitive, imprimitive
subgroup. So we can assume that $G$ is imprimitive. Then for some prime
divisor $p$ of $n$ there exists a partition $I=\{I_p^i:\, i\in [1,l]\}$ of
$V$ on $G$-isomorphic, automorphic $p$-subspaces. Let $X_n$ be a $n$-orbit
of $G$, $X_p=\hat{p}(I_p^i)X_n$ and $X_{2p}^{ij}=\hat{p}(I_p^i\circ
I_p^j)X_n=\uplus X_p\circ \uplus X_p$. Accordingly to corollary
\ref{mMmM->MM} $X_{2p}^{ij}=\cup_t (X_p\stackrel{\phi_t(i,j)}{\circ}
X_p)\equiv\cup L_p(i,j)$ for some maps $\{\phi_t(i,j)\}$. It follows that
the action of any cycle of length $p$, that is generated with some
$p$-rcycle from $X_p$, on partitions $L_p(i,j)$ generates cycles of
length $p$ that are again connected with $p$-rcycles from $X_p$. Hence $G$
contains a regular permutation of order $p$.


\section{Some applications of $k$-orbit theory}

\subsection{Solvability of groups of odd order}

Now we shall show that the $k$-orbit theory gives a simple proof of the
W. Feit, J.G. Thompson theorem: Solvability of groups of odd order
\cite{G}.

In this section we do not difference between a finite group and its
$ld$-representation. Also we assume that a finite group is not a direct
product.

\begin{lemma}\label{Aut(incoh)-not-simple}
Let $X_k\in Orb_k(G)$ be an incoherent $k$-orbit, then $G$ is not a
simple group.

\end{lemma}
\proof
Let $Y_k\subset X_k$ be a $k$-block and $L_k=Aut(X_k)Y_k$, then it is
evident that a concatenation of classes from $L_k$ is a $n$-orbit of a
normal subgroup $H\lhd Aut(X_k)=Aut(L_k)$ and that every transitive
subgroup $G<Aut(L_k)$ has non-trivial intersection with $H$. $\Box$

\begin{corollary}\label{impr-not.simple}
Let $G$ be a (non-cyclic) simple group, then it is (non-trivial) primitive.

\end{corollary}

\begin{theorem}\label{prim-invol.}
Any primitive group $G$ contains an involution.

\end{theorem}
\proof
Let $n$ be odd, then there exist odd numbers $k<l\leq n$, so that
$(k,l)=1$ and there exist automorphic subspaces $I_k\subset I_l$. If
$r=[l/k]$ is odd, then $m=l-rk$ is even and there exists an automorphic
subspace $I_m$ (lemma \ref{k-rorbit}). If $r$ is even, then $m$ is odd
and, because of primitivity of $G$, we can choose $l:=k$ and $k:=m$, if
$(k,m)=1$, or else $l:=l$ and $k:=m$. $\Box$

\begin{corollary}\label{odd-solv}
Let $G$ be a (non-cyclic) simple group, then it contains an involution.

\end{corollary}

\begin{corollary}\label{odd.od-impr}
Let $G$ be a group of odd order, then it is an imprimitive and hence
solvable group.

\end{corollary}

\subsection{A full invariant of a finite group}

Here we shall discuss the problem of a full invariant of a finite group
$F$ and assume that $F$ is not a direct product.

At first we can note that, if the ld-representation $G$ of $F$ is unique
accurate to similarity, then a full invariant of $F$ is defined by a full
invariant of $G$. Then we have two cases: $G$ is primitive and $G$ is
imprimitive.

\subsubsection{Let $G$ be primitive}

\begin{proposition}\label{k=n-l[n/l]}
Let $\neg(l|n)$, $I_l$ be an automorphic subspace and $k=n-l[n/l]$, then

\begin{enumerate}
\item\label{p^m<n}
$k$ divides $|G|$;

\item
there exists an automorphic $k$-subspace $I_k$;

\item
a subgroup $Stab(I_k)$ has non-trivial normalizer in $G$.

\end{enumerate}

\end{proposition}
\proof
The statements follows directly from lemmas \ref{k-rorbit} and
\ref{fix-tuple}. $\Box$

\begin{proposition}\label{V-coherent}
Let $\neg(k|n)$ and $X_k$ be a $k$-rorbit, then $X_k$ contains an
elementary $V$-coherent $k$-suborbit.

\end{proposition}
\proof
The statement follows from the definition of an elementary $V$-coherent
$k$-orbit. $\Box$\bigskip

So we see that $|G|$ and $n$ are high dependent and in general the degree
$n$ allows to define whether there can exist a group $G$ of order $\nu$.
Also an elementary $U$-coherent $k$-suborbit on every possible automorphic
subset $U\subset V$ is unique accurate to similarity.

All these facts suggest us the hypothesis that $|G|$ and $n$ could be a
full invariant of $G$ in the considered case.

\subsubsection{Let $G$ be imprimitive}

Let $k$ be a maximal automorphic divisor of $n$, then there exists an
incoherent $k$-rorbit $X_k$ of $G$. Let $Y_k\subset X_k$ be a $k$-block
and $L_k=Aut(X_k)Y_k=Aut(L_k)Y_k$, then $G<Aut(L_k)$, $L_k=GY_k$ and $Y_k$
is a $k$-orbit of a maximal normal subgroup $H\lhd G$, because of
lemma \ref{Aut(incoh)-not-simple}.

Let us to assume that we know $Y_k$, $|G|$ and $n$. It gives us the next
information: $|L_k|=n/k$, $|X_k|=|Y_k||L_k|$, $|H|=|G|/|L_k|$ and
$|Stab(I_k)|=|G|/|X_k|=|H|/|Y_k|$. In addition we know that $n$-orbit
$X_n$ of $G$ is a $|L_k|\times |L_k|$ matrix $M$, whose elements are
multi-classes of $L_k$ and which gives a regular representation of a
factor group $\Phi=G/H$. Since $H$ is a maximal normal subgroup, hence
$\Phi$ is a simple group.

Let $\Phi$ be not a cyclic group, then it is not a ld-representation. But
from here follows that $G$ is not a ld-representation too. This
contradiction shows us that $\Phi$ is always a cyclic group.

In order to give a full description of the group $G$ we have to find
the elements of matrix $M$. Let $L_k=\{Y_k^i:\,i\in [1,p]\}$ and
$\Phi=gr((1\ldots p))$, then we know that $M_{ij}=Y_k^r$, where
$r=i+j-1(mod\ p)$. One of possible construction of $M$ is obtained as
next. Every element $M_{ij}$ is obtained from the element $M_{1j}$ by
permutation of columns and every element $M_{1j}$ is obtained from the
element $M_{11}$ by permutation of lines. The columns of the element
$M{ij}$ are permutated relative to the columns of the element $M_{1j}$
with automorphisms of $M_{1j}$ that are similar to automorphisms of
$M_{11}$. The lines of the element $M_{1j}$ are permutated relative to
lines of the element $M_{11}$ on the condition to maintain the
automorphism property of $M_{11}$.

So, for obtaining of a full number invariant in this imprimitive case it
wants to find the number invariants that allow to calculate corresponding
$p^2-1$ permutations.

Now we consider the properties of $Y_k$.

\begin{lemma}\label{Y_k-prim}
$Aut(Y_k)$ is imprimitive.

\end{lemma}
\proof
Let $Aut(Y_k)$ be primitive, then $Y_k$ contains an elementary coherent
$q$-subrorbit for some prime $q<k$. From here follows that permutations
generating elements of matrix $M$ are trivial and hence $G$ is a direct
product. Contradiction. $\Box$\bigskip

\begin{corollary}\label{Y_k-reg}
$Aut(Y_k)$ is regular.

\end{corollary}

\begin{corollary}\label{p-group}
$G$ is $p$-group.

\end{corollary}

So we can formulate

\begin{hypothesis}\label{full.inv}
Full invariant of not $p$-group is defined with $|G|$ and $n$ and full
invariant of $p$-group of order $p^m$ is defined with maximum $mp^2$
permutations or corresponding numbers that define these permutations.

\end{hypothesis}

\subsection{The polynomial algorithm of graph isomorphism testing}

The graph isomorphism problem has a polynomial solution, if the problem of
separating of orbits of the automorphism group of a graph has a polynomial
solution. So we want to find the partition $O_2=Aut(X_2)X_2\subset
Orb_2(Aut(X_2))$ of a $2$-set $X_2\subset V^{(2)}$ polynomially on $n$.

Let $X_k\subset V^{(k)}$ be a $k$-set. We say that $X_k$ is
\emph{transitive}, if All $1$-projections of $X_k$, are equal. We say that
$X_k$ is \emph{regular}, if it satisfy to the two conditions:

\begin{enumerate}
\item
Every $l$-multiprojection of $X_k$ for $l\in [1,k]$ is homogeneous.

\item
All $l$-projections of $X_k$, containing the same $l$-tuple, are equal.

\end{enumerate}

\begin{lemma}\label{X_2=X_1}
Let $X_2$ be a regular $2$-set and $X_1^1,X_1^2$ be its $1$-projections.
Let $|X_2|=|X_1^1|$. If $X_1^1\neq X_1^2$, then $X_2$ is automorphic. If
$X_1^1=X_1^2$, then the partition $O_2=Aut(X_2)X_2$ is detected directly.

\end{lemma}

So the problem presents, if $|X_2|/|X_1^1|>1$ and $|X_2|/|X_1^2)|>1$. Let
$X_2\subset V^{(2)}$ be a regular $2$-set and $\cup Co(X_2)=V$, then
$Aut(X_2)$ is (generally intransitive) group of degree $n$ and all prime
cycles from $Aut(X_2)$ have the length not greater than $n$. Let $X_2$ be
automorphic and $Y_2\subset X_2$ be a $2$-orbit of a subgroup
$A<Aut(X_2)$, then $R_2=AX_2$ is a partition of $X_2$ on $2$-orbits of $A$
and hence classes of $R_2$ are regular $2$-sets.

The fundamental role in polynomial solution of considered problem plays the
theorem \ref{rr}. One can see that cyclic structure, that was described
for transitive $(2k)$-orbits, exists also on intransitive $(2k)$-orbits.
But the direction (left, right, left, right,\,\ldots)  must be change to
(left, right, right, left, left,\,\ldots). We can also note that in
general, if we have a regular $2$-set $X_2\subset V^{(2)}$, then we have
also a partition $P_2$ of $V^{(2)}$ on regular $2$-sets invariant relative
to $Aut(X_2)$. Thus with given intransitive regular $2$-set $X_2$ we can
also find transitive, regular, $Aut(X_2)$-invariant $2$-sets of $P_2$.

\begin{algorithm}\label{R_2,reg}
Let now $X_2$ be arbitrary regular, $2$-set and we try find its
automorphism, then we follow the next steps:

\begin{enumerate}
\item
Find an automorphic $2$-subset $Y_2\subset X_2$ that is expected to be a
$2$-suborbit of $Aut(X_2)$.

\item
Construct a partition $R_2^i$, $i=1,2,\ldots$, iterating the process that
follows from theorem \ref{rr}. By each iteration verify classes of $R_2^i$
on regularity and subpartition them if they are not regular.

\end{enumerate}
\end{algorithm}

This process leads to an automorphism, possibly trivial.

To find $Y_2$ is not difficult. At first it could be taken a subset
$Y_2(u)$ of $X_2$ whose $1$-projection $\hat{p}(I_1^1)Y_2$ is an element
$u$ of $V$.  Thus we can define whether $Stab(u)$ is trivial. If it is
trivial, then it follows that $X_2$ is incoherent and can be partition on
coherent $2$-subsets.

If $Aut(X_2)$ is trivial, then given algorithm detects this triviality in
maximally $n$ steps, if to assume that in each step only one element of
$V$ is separated.

Using lemma \ref{AV*BV,AV+BV} we can, having automorphic partitions $R_k'$
and $R_k''$ of $X_k$, obtain new more big partition $R_k=R_k'\sqcup
R_k''$.

For simplification of the process it can be chosen for partitioning on
$i$-th iteration the most suitable $2$-set $X_2^i$ from the partition
$P_2^i$ of $V^{(2)}$, on regular $2$-sets invariant to $Aut(X_2)$, and, by
partitioning of $X_2^i$, the whole partition $P_2^i$ can be further
partitioned to regular classes $P_2^{i+1}$ and used in the next iteration.


\section*{Conclusion}

This work was initiated by the polycirculant conjecture, described by P.
Comeron in his text \cite{Cameron} and represented on the site
(\url{http://www.maths.qmw.ac.uk/~pjc/homepage.html}).

The using of $k$-orbits to the polynomial solution of the graph isomorphism
problem was begun by Author in 1984. The generalization of $k$-orbits,
regular $k$-sets, was used for describing of the structure of strongly
regular graphs and their generalization on dimensions greater as two. This
approach discovered the difference between the structure of strongly
symmetrical but not automorphic partitions of $V^{(k)}$ and automorphic
partitions of $V^{(k)}$.

From this point of view the polycirculant conjecture seemed enough simple.
But nevertheless to find a correct proof was very difficult and only the
analysis of two examples of permutation groups (one elusive group of order
72 and degree 12, and the automorphism group of Peterson graph), that was
presented to author by P. Comeron, leaded to discovery of the specific
properties of $k$-orbits, not detectable with group theory, that brought a
proof of the conjecture.

In 1997 Author understood the connection between the graph isomorphism
problem and the problem of a full invariant of a finite group and has done
some attempts to obtain this full invariant by construction of some
appropriate group representations. This work gave better understanding of
the problem but did not bring the expected result. By construction of the
$k$-orbit theory it was of interest to consider a finite group with new
representation and this time the result was obtained.

Also the specific symmetry properties of $k$-orbits, that are not visible
in other most algebraic theories, gave possibility for simple polynomial
solution of the graph isomorphism problem.


\section*{Acknowledgements}
I would like to express many thanks to Prof.~P.~Cameron for encouraging to
do this work and for the counterexamples to not correct hypotheses on the
way to prove the polycirculant conjecture.



\begin{thebibliography}{9}
\bibitem{Cameron} Peter J. Cameron, \emph{Elusive groups and the
polycirculant conjecture}, Queen Mary, University of London, Barcelona,
February 2001

\bibitem{G} Daniel Gorenstein, \emph{Finite simple groups}, Plenum Press,
New York and London, 1982

\end{thebibliography}
\end{document}